\documentclass[10pt, a4paper]{article}
\usepackage{comment}
\usepackage{graphics}
\usepackage{latexsym}
\usepackage{amssymb}
\usepackage{amsmath}
\usepackage{color}
\usepackage{epsfig}
\usepackage{pgf,tikz}
\usepackage{longtable}
\newcommand{\beq}{\begin{equation}}
\newcommand{\eeq}{\end{equation}}
\newcommand{\beqas}{\begin{eqnarray*}}
\newcommand{\eeqas}{\end{eqnarray*}}
\newcommand{\ep}{\varepsilon}

\newcommand{\ue}{u^{\varepsilon}}

\newcommand{\bY}{{\bf Y}}
\newcommand{\by}{\boldsymbol{y}}

\newcommand{\bV}{{\bf V}}
\newcommand{\bbV}{{\bar{\bf V}}}
\newcommand{\hbV}{{\hat{\bf V}}}

\newcommand{\bv}{\boldsymbol{v}}
\newcommand{\bu}{\boldsymbol{u}}

\newcommand{\bbu}{\bar{\boldsymbol{u}}}
\newcommand{\hbu}{\hat{\boldsymbol{u}}}
\newcommand{\bbv}{\bar{\boldsymbol{v}}}
\newcommand{\hbv}{\hat{\boldsymbol{v}}}

\newcommand{\bau}{\bar u}
\newcommand{\hau}{\hat u}
\newcommand{\bav}{\bar v}
\newcommand{\hav}{\hat v}
\newcommand{\bafu}{\bar{\frak u}}
\newcommand{\hafu}{\hat{\frak u}}
\newcommand{\bafv}{\bar{\frak v}}
\newcommand{\hafv}{\hat{\frak v}}

\newcommand{\bw}{\boldsymbol{w}}

\newcommand{\bH}{\boldsymbol{\cal H}}

\newcommand{\boldp}{\boldsymbol{p}}
\newcommand{\boldq}{\boldsymbol{q}}

\newcommand{\tL}{\tilde{L}^\infty}

\newcommand{\fu}{\frak u}
\newcommand{\fv}{\frak v}
\newcommand{\fw}{\frak w}
\newcommand{\fp}{\frak p}
\newcommand{\fq}{\frak q}

\newcommand{\bt}{\begin{theorem}}
\newcommand{\bp}{\begin{proposition}}
\newcommand{\bl}{\begin{lemma}}
\newcommand{\br}{\begin{remark}}

\newcommand{\et}{\end{theorem}}
\newcommand{\epr}{\end{proposition}}
\newcommand{\el}{\end{lemma}}
\newcommand{\er}{\end{remark}}

\def\curl{{\rm curl\,}}
\def\scurl{{\rm curl}}

\newcommand{\bproof}{{\it Proof}\ \ }
\newcommand{\eproof}{{\hfill$\Box$}}

\newcommand{\be}{\begin{equation}}
\newcommand{\ee}{\end{equation}}
\newcommand{\bea}{\begin{eqnarray}}
\newcommand{\eea}{\end{eqnarray}}
\newcommand{\beas}{\begin{eqnarray*}}
\newcommand{\eeas}{\end{eqnarray*}}

\def\IR{\mathbb{R}}

\def\curl{{\rm curl\,}}

\def\scurl{{\rm curl}}

\newtheorem{theorem}{Theorem}[section]
\newtheorem{lemma}[theorem]{Lemma}
\newtheorem{assumption}[theorem]{Assumption}
\newtheorem{proposition}[theorem]{Proposition}
\newtheorem{definition}[theorem]{Definition}
\newtheorem{remark}[theorem]{Remark}

\definecolor{darkgreen}{rgb}{0,0.7,0}

\numberwithin{equation}{section}
\numberwithin{equation}{section}
\title{High dimensional finite elements for multiscale Maxwell wave equations}
\author{
Van Tiep Chu
and 
Viet Ha Hoang\\[10pt]
Division of Mathematical Sciences,\\
School of Physical and Mathematical Sciences,\\
 Nanyang Technological University, Singapore 637371 
}

\date{}
\begin{document}
\maketitle

\begin{abstract}
We develop an essentially optimal numerical method for solving multiscale Maxwell wave equations in a domain $D\subset\IR^d$. The problems depend on $n+1$ scales: one macroscopic scale and $n$ microscopic scales. Solving the macroscopic multiscale homogenized problem, we obtain the desired macroscopic and microscopic information. This problem depends on $n+1$ variables in $\IR^d$, one for each scale that the original multiscale equation depends on, and is thus posed in a high dimensional tensorized domain. The straightforward full tensor product finite element (FE) method is exceedingly expensive. We develop the sparse tensor product FEs that solve this multiscale homogenized problem with essentially optimal number of degrees of freedom, that is essentially equal to that required for solving a macroscopic problem in a domain in $\IR^d$ only, for obtaining a required level of accuracy. Numerical correctors are constructed from the FE solution. For two scale problems, we derive a rate of convergence for the numerical corrector in terms of the microscopic scale and the FE mesh width. Numerical examples confirm our analysis.  
\end{abstract}

\section{Introduction}
We study the high dimensional finite element (FE) method for solving multiscale Maxwell wave equations in a domain $D\subset\IR^d$. The equation depends on the macroscopic scale and $n$ microscopic scales, and is locally periodic. We study the problem via multiscale convergence. In the limit where all the microscopic scales converge to zero, we obtain the multiscale homogenized equation. This equation contains the solution to the homogenized equation which approximates the solution of the original multiscale equation macroscopically, and the scale interacting corrector terms which provide the microscopic behaviour of the solution. Solving the equation, we obtain all the necessary  information. However, the multiscale homogenized equation is posed in a high dimensional tensorized domain. It depends on $n+1$ variables in $\IR^d$, one for each scale. The direct full tensor product FE method is highly expensive. We develop the sparse tensor product FEs to solve this problem which requires only essentially equal number of degrees of freedom as for solving an equation posed in $\IR^d$ for obtaining a required level of accuracy. The complexity is thus essentially optimal. 

As for any other multiscale problems, a direct numerical method using fine mesh to capture the microscopic scales is prohibitively expensive. There have been attempts to develop numerical methods for solving multiscale wave equations, and multiscale Maxwell equations with reduced complexity, though comparing to other types of multiscale equations, multiscale wave and multiscale Maxwell equations have been paid far less attention.

For multiscale wave equations, in \cite{Owhadi2008} Owhadi and Zhang  build a set of basis functions that contain microscopic information from the solutions of $d$ multiscale equations. These equations are solved using fine mesh to capture the microscopic scales. In \cite{JiangEfendievGinting}, Jiang et al. employ the Multiscale Finite Element method (\cite{HouWu}, \cite{EfendievHou}) to solve wave equations that depend on a continuum spectrum of scales, using limited global information. The Heterogeneous Multiscale Method (HMM) (\cite{EEngquist}, \cite{Abdulleetalacta}) is employed by Engquist et al. using finite differences to solve multiscale wave equations that show the dispersive behaviour at large time. Abdulle and Grote \cite{AbdulleGrote} employ the Heterogeneous Multiscale Method (HMM) to solve multiscale equation using finite elements. The approaches in these papers are general, but the complexity at each time step grows superlinearly with respect to the optimal complexity level. In \cite{XHwave}, Xia and Hoang develop the essentially optimal sparse tensor product FE method for locally periodic multiscale wave equations; the complexity of the method only grows log-linearly at each time step. The method is employed successfully for multiscale elastic wave equations in \cite{XHelasticwave}.

There has not been much research on efficient numerical methods for multiscale Maxwell equations. The traditional method that constructs the homogenized equation by solving cell problems is considered in \cite{Zhangetal2017} (see also the related references therein) where a set of cell problems are solved at each macroscopic points. The complexity is thus very high. The HMM method is applied for multiscale Maxwell equations in frequency domain in Ciarlet et al. \cite{CiarletMaxwell}. Ohlberger et al. considered a locally periodic two sale harmonic Maxwell equation in \cite{Ohlbergermaxwell}  though the problem is assumed uniformly coercive with respect to the microscopic scale. The HMM method is analyzed for the two scale homogenized problem using the approach in \cite{Ohlberger}. The complexity of the method is equivalent to that of a full tensor product FE method for solving the two scale homogenized equation. In \cite{Tiep1}, Chu and Hoang develop the sparse tensor product edge FE method for locally periodic stationary multiscale Maxwell equations. The method requires only a number of degrees of freedom that is essentially equivalent to that needed for solving a macroscopic scale Maxwell equation in a domain in  $\IR^d$, and is therefore optimal. Chu and Hoang \cite{Tiep1} construct numerical correctors from the finite element solutions. For two scale problems, an explicit error in terms of the FE error and the homogenization error is deduced for the numerical corrector.

%For The sparse tensor product finite element has been successfully applied for other types of multiscale equations in \cite{HSelliptic}, \cite{Hoangmonotone}, \cite{XHelasticity}. 
We develop the sparse tensor product FE approach for multiscale Maxwell wave equations in this paper using edge FEs. We show that the complexity of the method is essentially optimal. The sparse tensor product FE approach for multiscale problems is initiated by Hoang and Schwab in \cite{HSelliptic} for elliptic equations, and is applied for other types of equations in \cite{Hoangmonotone}, \cite{XHwave}, \cite{XHelasticity}, \cite{XHelasticwave}. 

In the next section, we set up the multiscale Maxwell wave equation and derive the multiscale homogenized equation. We will only summarize the results and refer to \cite{Mwt} for detailed derivation. In Section 3, we study FE approximation for the multiscale homogenized Maxwell wave equation using general FE spaces. In Subsection 3.1, we study the spatially semidiscrete problem where only the spatial variable is discretized. We follow the framework of Dupont \cite{Dupont} for wave equations. 
%For scalar multiscale wave equations, t
The approach has been applied for the multiscale homogenized equations of scalar multiscale wave equations in Xia and Hoang \cite{XHwave}. However, the application of the framework to multiscale homogenized Maxwell wave equations requires substantial modification for the analysis of the convergence due to the corrector terms $\fu_i$ in \eqref{eq:l1}. In Subsection 3.2, we consider the fully discrete problem where both the temporal and spatial variables are discretized. The convergence of the general discretization schemes in Section 3, and the full and sparse tensor product FE approximations in Section 5 require regularity for the solution of the multiscale homogenized Maxwell wave equation. In Section 4, we prove that the required regularity hold under mild conditions. In Section 5, we apply the discretization schemes in Section 3 for the full tensor product and the sparse tensor product edge FEs. We prove that the sparse tensor product FE method obtains an approximation with essentially the same level of accuracy as the full tensor product FEs but requires only essentially the same number of degrees of freedom as for solving a macroscopic Maxwell equation in $\IR^d$, and is thus essentially optimal. In Section 6, we construct numerical correctors from the FE solutions. For two scale problems, an explicit homogenization error in terms of the microscopic scale is available. The derivation is complicated, especially due to the low regularity of the solution of the homogenized Maxwell wave equation. We therefore only summarize the theoretical results  and refer the reader to \cite{Mwt} for details. From this, we derive a numerical corrector with an explicit error in terms of the homogenization error, and the FE error. For general multiscale problems, such a homogenization error is not available. We thus derive a general numerical corrector without an explicit error. Section 7 presents some numerical examples in two dimensions that confirm our analysis.  

Throughout the paper, by $\curl$ and $\nabla$ without explicitly indicating the variable, we mean the $\curl$ and the gradient of a function of $x$ with respect to $x$, and by $\curl_x$ and $\nabla_x$ we denote the partial $\curl$ and partial gradient of a function that depends on $x$ and other variables. We denote by $\langle\cdot,\cdot\rangle_{X',X}$ the duality pairing of a Banach space $X$ and its dual $X'$. Repeated indices indicate summation. The notation $\#$ denotes spaces of periodic functions with the period being the unit cube in $\IR^d$. 

\section{Multiscale Maxwell wave problems}
We set up the multiscale Maxwell wave equation and use multiscale convergence to homogenize it in this section. 
\subsection{Problem setting}
Let $D$ be a bounded domain in $\IR^d$ ($d=2,3$). Let $Y$ be the unit cube in $\IR^d$. By $Y_1,\ldots,Y_n$ we denote $n$ copies of $Y$. We denote by $\bY$ the product set $Y_1\times Y_2\times\ldots\times Y_n$ and by $\by=(y_1,\ldots,y_n)$. For $i=1,\ldots,n$, we denote by $\bY_i=Y_1\times\ldots\times Y_i$. Let $a$ and $b$ be  functions from $D\times Y_1\times\ldots\times Y_n$ to $\IR^{d\times d}_{sym}$. We assume that $a$ and $b$ satisfy the boundedness and coerciveness conditions: for all $x\in D$ and $\by\in\bY$, and all $\xi,\zeta\in \IR^d$,
\be
\begin{array}{lr}
\displaystyle{\alpha|\xi|^2\le a_{ij}(x,\by)\xi_i\xi_j,\ \ \ a_{ij}\xi_i\zeta_j\le\beta|\xi||\zeta|}\\
\displaystyle{\alpha|\xi|^2\le b_{ij}(x,\by)\xi_i\xi_j,\ \ \ b_{ij}\xi_i\zeta_j\le\beta|\xi||\zeta|}
\end{array}
\label{eq:coercive}
\ee
where $\alpha$ and $\beta$ are positive numbers. 
 Let $\ep$ be a small positive value, and $\ep_1,\ldots,\ep_n$ be $n$ functions of $\ep$ that denote the $n$ microscopic scales that the problem depends on. We assume the following scale separation properties: for all $i=1,\ldots,n-1$
\be
\lim_{\ep\to 0}{\ep_{i+1}(\ep)\over\ep_i(\ep)}=0.
\label{eq:scaleseparation}
\ee
Without loss of generality, we assume that $\ep_1(\ep)=\ep$. 
We define the multiscale coefficients of the Maxwell equation $a^\ep$ and $b^\ep$ which are functions from $D$ to $\IR^{d\times d}_{sym}$ as 
\[
a^\ep(x)=a(x,{x\over\ep_1},\ldots,{x\over\ep_n}),\ \ b^\ep(x)=b(x,{x\over\ep_1},\ldots,{x\over\ep_n}).
\]
When $d=3$ we define the space
\[
W=H_0(\curl,D)=\{u\in L^2(\Omega)^3,\ \ \curl u\in L^2(\Omega)^3,\ \ u\times\nu=0\},\ \ 
H=L^2(D)^3
\]
and when $d=2$
\[
W=H_0(\curl,D)=\{u\in L^2(D)^2,\ \ \curl u\in L^2(D),\ \ u\times\nu=0\},\ \ 
H=L^2(D)^2
\]
where $\nu$ denotes the outward normal vector on the boundary $\partial D$. These spaces form the Gelfand triple $W\subset H\subset W'$. 
%We denote by $\langle\cdot,\cdot\rangle$ the inner product in $H$, extending to the duality pairing between $W'$ and $W$. 
 We note that when $d=3$, $\curl\ue$ is a vector function in $L^2(D)^3$ and when $d=2$, $\curl\ue$ is a scalar function in $L^2(D)$. 
Let $f\in L^2(0,T; H)$, $g_0\in W$ and $g_1\in H$. 
We consider the problem:
Find $u^\varepsilon (t,x) \in L^2(0,T; W) $ so that
\[
\begin{cases}
b^\varepsilon(x)\dfrac{\partial^2u^\varepsilon(t,x) }{\partial t^2}+ \mathrm{curl}( a^\varepsilon(x)\mathrm{curl} u^\varepsilon(t,x) )=f(t,x),&  (0,T) \times D\\
%\nabla \cdot u^\varepsilon =0,&(0,T) \times D\\
%u \times n=0& \\
u^\varepsilon(0,x)=g_0(x)&\\
u_t^\varepsilon(0,x)=g_1(x)&\\

\end{cases}
\]
with the boundary condition $\ue\times\nu=0$ on $\partial D$.  We will mostly present the analysis for the case $d=3$ and only discuss the case $d=2$ when there is significant difference.  For notational conciseness, we denote by
\be
H_i=L^2(D\times \bY_i)^3,\ \ i=1,\ldots,n.
\label{eq:H}
\ee
In variational form, this problem becomes: Find $\ue\in L^2(0,T; W) \cap H^1(0,T;H)$ so that
\be
\left\langle b^\varepsilon(x)\frac{\partial^2 \ue}{\partial t^2}, \phi(x)\right\rangle_{W',W}+\int_Da^\varepsilon(x)\mathrm{curl} u^\varepsilon(t,x)\cdot\curl\phi(x)dx
=\int_Df(t,x)\cdot\phi(x)dx
\label{eq:varprob}
\ee
for all $\phi\in W$ when $d=3$; and when $d=2$ we need to replace the vector product for $\curl$ by the scalar multiplication. 
Problem \eqref{eq:varprob} has a unique solution $\ue\in L^2(0,T;W)\bigcap H^1(0,T;H)\bigcap H^2(0,T;W')$ that satisfies
\be
\|\ue\|_{L^2(0,T;W)}+\|\ue\|_{H^1(0,T;H)}+\|\ue\|_{H^2(0,T;W')}\le c(\|f\|_{L^2(0,T;H)}+\|g_0\|_W+\|g_1\|_H)
\label{eq:uepbound}
\ee
where the constant $c$ only depends on the constants $\alpha$ and $\beta$ in \eqref{eq:coercive} and $T$ (see Wloka \cite{Wloka}). 

We will study this problem via multiscale convergence.
%
%\subsection{Multiscale covergence}
%
\subsection{Multiscale convergence}
We study homogenization of problem \eqref{eq:varprob} via multiscale convergence. We therefore recall the definition of multiscale convergence (see Nguetseng \cite{Nguetseng1989}, Allaire \cite{Allaire1992} and Allaire and Briane \cite{Allaire1996}). 
\begin{definition}
A sequence of functions $\{w^\ep\}_\ep\subset L^2(0,T;H)$ $(n+1)$-scale converges to a function $w_0\in L^2(0,T;D\times \bY)$ if for all smooth functions $\phi(t,x, \by)$ which are $Y$ periodic w.r.t $y_i$ for all $i=1,\ldots,n$:
\[
\lim_{\ep \to 0}\int_0^T\int_Dw^\ep(t,x)\phi(t,x,{x\over\ep_1},\ldots,{x\over\ep_n})dxdt=\int_0^T\int_D\int_\bY w_0(t,x,\by)\phi(t,x,\by)d\by dxdt.
\]
\end{definition}
We have the following result.
\begin{proposition}\label{prop:2scexistence}
From a bounded sequence in $L^2(0,T; H)$ we can extract an $(n+1)$-scale convergent subsequence.
\end{proposition}
We note that the definition above for functions which depend also on $t$ is slightly different from that in \cite{Nguetseng1989} and \cite{Allaire1992} as we take also the integral with respect to $t$. However, the proof of Proposition \ref{prop:2scexistence} is similar.

For a bounded sequence in $L^2(0,T;W)$, we have the following results  which are very similar to those in \cite{Tiep1} and \cite{Wellander} for functions which do not depend on $t$. The proofs for these results are very similar to those in \cite{Tiep1} so we do not present them here. As in \cite{Tiep1}, we denote by $\tilde H_\#(\curl,Y_i)$ the space of equivalent classes of functions in $H_\#(\curl,Y_i)$ of equal $\curl$. 
\begin{proposition}\label{prop:msHcurl}
Let $\{w^\ep\}_\ep$ be a bounded sequence in  $L^2(0,T; W)$. There is a subsequence (not renumbered), a function $w_0\in  L^2(0,T; W)$, $n$ functions $\frak w_i\in L^2((0,T) \times D\times Y_1\times\ldots\times Y_{i-1},H^1_\#(Y_i)/\IR)$ such that 
\[
w^\ep \stackrel{(n+1)-\mbox{scale}}{\longrightarrow}w_0+\sum_{i=1}^n\nabla_{y_i}\frak w_i.
\]
Further, there are $n$ functions $w_i\in L^2((0,T) \times D\times\ldots\times Y_{i-1},\tilde H_\#(\curl,Y_i))$ such that
\[
\curl w^\ep \stackrel{(n+1)-\mbox{scale}}{\longrightarrow}\curl w_0+\sum_{i=1}^n\scurl_{y_i}w_i.
\]
\end{proposition}
From \eqref{eq:uepbound} and Proposition \eqref{prop:msHcurl}, we can extract a subsequence (not renumbered), a function $u_0\in L^2(0,T; W)$, $n$ functions $\frak u_i\in L^2(0,T;D\times Y_1\times\ldots\times Y_{i-1},H^1_\#(Y_i)/\IR)$ and $n$ functions $u_i\in L^2(0,T;D\times Y_1\times\ldots\times Y_{i-1},\tilde H_\#(\curl,Y_i))$ such that
\be
\ue\stackrel{(n+1)-\mbox{scale}}{\longrightarrow} u_0+\sum_{i=1}^n\nabla_{y_i}\frak u_i,
\label{eq:l1}
\ee
and 
\be
\curl\ue \stackrel{(n+1)-\mbox{scale}}{\longrightarrow} \curl u_0+\sum_{i=1}^n\scurl_{y_i}u_i.
\label{eq:l2}
\ee
For $i=1,\ldots,n$, let $W_i=L^2(D\times Y_1\times\ldots\times Y_{i-1},\tilde H_\#(\curl,Y_i))$ and $V_i=L^2(D\times Y_1\times\ldots\times Y_{i-1},H^1_\#(Y_i)/\IR)$. We define the space $\bV$ as
\be
\bV=W\times W_1\times\ldots\times W_n\times V_1\times\ldots\times V_n.
\label{eq:spacebV}
\ee
For $\bv=(v_0,\{v_i\},\{\frak v_i\})\in \bV$, we define the norm
\[
|||\bv|||=\|v_0\|_{H(\curl,D)}+\sum_{i=1}^n\|v_i\|_{L^2(D\times\bY_{i-1},\tilde H_\#(\curl,Y_i))}+\sum_{i=1}^n\|\frak v_i\|_{L^2(D\times\bY_{i-1},H^1_\#(Y_i))}.
\]
Let $\bu=(u_0,\{u_i\}, \{\frak u_i\}) \in \bV$. We define the function
\[
\int_{\bY}b(x,\by)\left(u_0(t,x)+\sum_{i=1}^n\nabla_{y_i}\fu_i(t,x,\by_i)\right)d\by
\]
in $W'$ so that
\[
\left\langle\int_{\bY}b(x,\by)\left(u_0(t,x)+\sum_{i=1}^n\nabla_{y_i}\fu_i(t,x,\by_i)\right)d\by,v_0\right\rangle_{W',W}=\int_D\int_\bY b(x,\by)\left(u_0(t,x)+\sum_{i=1}^n\nabla_{y_i}\fu_i(t,x,\by_i)\right)\cdot v_0d\by dx;
\]
and the function $b(x,\by)(u_0(t,x)+\sum_{i=1}^n\nabla_{y_i}\fu_i(t,x,\by))$ in $W_j'$ as
\[
\left\langle b(x,\by)\left(u_0(t,x)+\sum_{i=1}^n\nabla_{y_i}\fu_i(t,x,\by_i)\right),\fv_j\right\rangle_{V_j',V_j}=\int_D\int_\bY b(x,\by)(u_0(t,x)+\sum_{i=1}^n\nabla_{y_i}\fu_i(t,x,\by_i))\cdot\nabla_{y_j}\fv_j d\by dx.
\]
We then have the following result.
\begin{proposition}\label{prop:1}
The function $\bu=(u_0,\{u_i\},\{\bu_i\})$ satisfies
%\beqas
%&&\left({\partial^2\over\partial t^2}\int_{\bY}b(x,\by)\left(u_0(t,x)+\sum_{i=1}^n\nabla_{y_i}\fu_i(t,x,\by)\right)d\by,v_0\right)_{W',W}=\\
%&&\int_D f(x)\cdot v_0(x)dx-\int_D\int_\bY a(x,\by)\left(\curl u_0+\sum_{i=1}^n\scurl_{y_i}u_i\right)\cdot\left( v_0+\sum_{i=1}^n\scurl_{y_i}v_i\right)d\by dx
%\eeqas
%and 
%\[
%\left({\partial^2\over\partial t^2}\left(b(x,\by)\left(u_0(t,x)+\sum_{i=1}^n\nabla_{y_i}\fu_i(t,x,\by)\right)\right),\fv_i\right)_{W_i',W_i}=0,
%\]
%i.e. the function $\bu$ satisfies the multiscale homogenized equation
\begin{multline}
\left\langle{\partial^2\over\partial t^2}\int_\bY b(x,\by)\left(u_0(t,x)+\sum_{i=1}^n\nabla_{y_i}\fu_i(t,x,\by_i)\right)d\by,v_0\right\rangle_{W',W}+\\
\sum_{j=1}^n\left\langle{\partial^2\over\partial t^2}b(x,\by)\left(u_0(t,x)+\sum_{i=1}^n\nabla_{y_i} \frak u_i(t,x,\by_i)\right),\frak v_j\right\rangle_{V_j',V_j}+\\
\int_D\int_{\bY}a(x,\by)\left(\curl u_0+\sum_{i=1}^n\scurl_{y_i}u_i\right)\cdot\left(\curl v_0+\sum_{i=1}^n\scurl_{y_i}v_i\right)d\by dx=\int_D f(t,x)\cdot v_0(x)dx
\label{eq:msprob}
\end{multline}
for all $\bv=(v_0,\{v_i\},\{\frak v_i\})\in \bV$. 
\end{proposition} 
For the initial conditions, we have
%%%%%%%%%%%%%%%%
\bp\label{prop:2} We have $u_0\in H^1(0,T;H)$, $\nabla_{y_i}\frak u_i\in H^1(0,T;H_i)$ for all $i=1,\ldots,n$. Further
\be
u_0(0,\cdot)=g_0,\ \ \nabla_{y_i}\frak u_i(0,\cdot,\cdot)=0,
\label{eq:u0i}
\ee
\be
{\partial\over\partial t}\int_\bY b(x,\by)(u_0(t,x)+\sum_{i=1}^n\nabla_{y_i}\fu_i(t,x,\by_i))d\by\bigg|_{t=0}=\int_\bY b(x,\by)g_1(x)d\by,\ \ \mbox{in}\ W'
\label{eq:fuiv0i}
\ee
and for $j=1,\ldots,n$
\be
{\partial\over\partial t}b(x,\by)(u_0(t,x)+\sum_{i=1}^n\nabla_{y_i}\fu_i(t,x,\by_i))\bigg|_{t=0}=b(x,\by)g_1(x),\ \ \mbox{in}\ V_j'.
\label{eq:fuivii}
\ee
\epr
We then have
\bp\label{prop:3} With the initial conditions \eqref{eq:u0i}, \eqref{eq:fuiv0i} and \eqref{eq:fuivii}, problem \eqref{eq:msprob} has a unique solution.
\epr

%For $\bv=(v_0,v_1,\ldots,v_n,\fv_1,\ldots,\fv_n)$ and $\bw=(w_0,w_1,\ldots,w_n,\fw_1,\ldots,\fw_n)$ in $\bV=W\times W_1\times\ldots W_n\times V_1,\ldots,V_n$. We define the bilinear forms
%\[
%A(\bv,\bw)=\int_D\int_\bY a(x,\by)(\curl v_0+\sum_{i=1}^n\scurl_{y_i}v_i)\cdot(\curl w_0+\sum_{i=1}^n\scurl_{y_i}w_i)d\by dx,
%\]
%and
%\[
%B(\bv,\bw)=\int_D\int_\bY b(x,\by)(v_0+\sum_{i=1}^n\nabla_{y_i}\fv_i)\cdot (w_0+\sum_{i=1}^n\nabla_{y_i}\fw_i)d\by dx.
%\]
%
The proofs of Propositions \ref{prop:1}, \ref{prop:2} and \ref{prop:3} can be found in \cite{Mwt}.
\section{Finite element discretization}
We study finite element approximation for problem \eqref{eq:msprob} in this section. We first consider the semidiscrete problem where we discretize the spatial variables. We then consider the fully discrete problem where both the temporal and spatial variables are discretized. 
\subsection{Spatially semidiscrete problem}
We consider in this section the spatial semidiscretization of the homogenized problem \eqref{eq:msprob}. For approximating $u_0$, we suppose that there is a hierarchy of finite dimensional subspaces 
\[
W^1\subset W^2\subset\ldots\subset W^L\ldots\subset W;
\]
to approximate $u_i$, $i=1,2,\ldots,n$, we assume a hierarchy of finite dimensional subspaces
\[
W_i^1\subset W_i^2\subset\ldots\subset W_i^L\ldots\subset W_i;
\]
and to approximate $\fu_i$, $i=1,2,\ldots,n$, we assume a hierarchy of finite dimensional subspaces
\[
V_i^1\subset V_i^2\subset\ldots\subset V_i^L\ldots\subset V_i.
\]
Let
\[
\bV^L=W^L\times W_1^L\times\ldots\times W_n^L\times V_1^L\times\ldots\times V_n^L
\]
which is a finite dimensional subspace of $\bV$ defined in \eqref{eq:spacebV}. 
We consider the spatially semidiscrete approximating problems: Find $\bu^L(t)=(u_0^L, u_1^L,...,u_n^L,\frak u_1^L,...,\frak u_n^L) \in \bV^L$ so that
\begin{multline}
\int_D\int_\bY \left[ b(x,\by)\left(\frac{\partial ^2}{\partial t^2}u^L_0(t,x)+\sum_{i=1}^n\nabla_{y_i}\frac{\partial ^2}{\partial t^2} \frak u^L_i(t,x,\by_i)\right)\cdot\left(v_0^L+\sum_{i=1}^n\nabla_{y_i} \frak v^L_i\right)\right.\\
\left.+a(x,\by)\left(\curl u^L_0+\sum_{i=1}^n\scurl_{y_i}u^L_i\right)\cdot\left(\curl v^L_0+\sum_{i=1}^n\scurl_{y_i}v^L_i\right)\right]d\by dx\\=\int_D f(t,x)\cdot v^L_0(x)dx
\label{eq:nmsprob}
\end{multline}
for all $\bv^L=(v_0^L,v_1^L,\ldots,v_n^L,\fv_0^L,\ldots,\fv_n^L) \in \bV^L$. Let $g_0^L\in W^L$, $g^L_1\in W^L$ which are approximations of $g_0$ and $g_1$ in $W$ and in $H$ respectively. The initial conditions \eqref{eq:u0i} are approximated by:
\be
u_0^L(0,\cdot)=g_0^L,\ \ \nabla_{y_i}\fu^L_i(0,\cdot,\cdot)=0.
\label{eq:u0iL}
\ee
We approximate the initial conditions \eqref{eq:fuiv0i} and \eqref{eq:fuivii} by
\begin{multline*}
\int_D\int_\bY b(x,\by)\left({\partial u_0^L\over\partial t}(0)+\sum_{i=1}^n{\partial\over\partial t}\nabla_{y_i}\fu_i^L(0)\right)\cdot\left(v_0^L+\sum_{i=1}^n\nabla_{y_i}\fv_i^L\right)d\by dx\\=\int_D\int_\bY b(x,\by)g_1^L(x)\cdot \left(v_0^L+\sum_{i=1}^n\nabla_{y_i}\fv_i^L\right)d\by dx
\end{multline*}
for all $v_0^L\in W^L$ and $\fv^L_i\in V_i^L$, i.e., 
\[
\int_D\int_\bY b(x,\by)\left({\partial u_0^L\over\partial t}(0)-g_1^L+\sum_{i=1}^n{\partial\over\partial t}\nabla_{y_i}\fu_i^L(0)\right)\cdot\left(v_0^L+\sum_{i=1}^n\nabla_{y_i}\fv_i^L\right)d\by dx=0.
\]
Using the coercivity of the matrix $b(x,\by)$, we get 
\be
{\partial u_0^L\over\partial t}(0)=g_1^L,\ \  {\partial\over\partial t}\nabla_{y_i}\fu_i^L(0)=0.
\label{eq:fu0fuiL}
\ee
For $\bv=(v_0,v_1,\ldots,v_n,\fv_1,\ldots,\fv_n)$ and $\bw=(w_0,w_1,\ldots,w_n,\fw_1,\ldots,\fw_n)$ in $\bV=W\times W_1\times\ldots W_n\times V_1\times\ldots\times V_n$, we define the bilinear forms
\[
A(\bv,\bw)=\int_D\int_\bY a(x,\by)\left(\curl v_0+\sum_{i=1}^n\scurl_{y_i}v_i\right)\cdot\left(\curl w_0+\sum_{i=1}^n\scurl_{y_i}w_i\right)d\by dx,
\]
and
\[
B(\bv,\bw)=\int_D\int_\bY b(x,\by)\left(v_0+\sum_{i=1}^n\nabla_{y_i}\fv_i\right)\cdot \left(w_0+\sum_{i=1}^n\nabla_{y_i}\fw_i\right)d\by dx.
\]
\begin{proposition}
Problem \eqref{eq:nmsprob} together with the initial conditions \eqref{eq:u0iL} and \eqref{eq:fu0fuiL} has a unique solution.
\end{proposition}
\bproof
In the bilinear form $B$, let $R$ be the  matrix that describes the interaction of the basis functions of $W^L$ with themselves, let $N$ be the  matrix that describes the interaction of the basis functions of $V_1^L\times\ldots\times V_n^L$ with themselves, and let $S$ be the matrix that describes the interaction of the basis functions of $W^L$ and the basis functions of $V_1^L\times\ldots\times V_n^L$. For the bilinear form $A$, let $Q$ be the matrix that describes the interaction of the basis functions of $W^L$ with themselves, let $P$ be the matrix describing the interaction of the basis functions of $W^L$ and $W_1^L\times\ldots\times W_n^L$, and let $M$ be the matrix describing the interactions of the basis functions of $W_1^L\times\ldots\times W_n^L$ and themselves. Let $F$ be the column vector describing the interaction of $f$ and the basis functions of $W^L$. Let $C_0$ be the coefficient vector in the expansion of $u_0^L$ with respect to the basis functions of $W^L$. Let $C_1$ be the coefficient vector in the expansion of $(u_1^L,\ldots,u_n^L)$ with respect to the basis functions of $W_1^L\times\ldots\times W_n^L$. Let ${\frak C}_1$ be the coefficient vector in the expansion of $(\fu_1^L,\ldots,\fu_n^L)$ with respect to the basis functions of $V_1^L\times\ldots\times V_n^L$. We have the following equations
\beqas
R{d^2C_0\over dt^2}+S{d^2{\frak C}_1\over dt^2}+QC_0+PC_1&=&F,\\
P^{\top}C_0+MC_1&=&0,\\
S^{\top}{d^2C_0\over dt^2}+N{d^2{\frak C}_1\over dt^2}&=&0.
\eeqas
Using  $C_1=-M^{-1}P^{\top}C_0$, we deduce the system
\beqas
\begin{bmatrix}
R & S\\
S^{\top} & N
\end{bmatrix}
{d^2\over dt^2}
\begin{bmatrix}
C_0\\
{\frak C}_1
\end{bmatrix}
+
\begin{bmatrix}
Q-PM^{-1}P^{\top} & 0\\
0 & 0
\end{bmatrix}
\begin{bmatrix}
C_0\\
{\frak C}_1
\end{bmatrix}
=
\begin{bmatrix}
F\\
0
\end{bmatrix}.
\eeqas
We note that  $\begin{bmatrix} R & S\\ S^{\top} & N \end{bmatrix}$ is the Gram  matrix for the interaction of the basis of $W^L$ and $V_1^L\times\ldots\times V_n^L$ in the bilinear form $B$ so is positive definite. The system thus has a unique solution.  
\eproof

For each $t \in (0,T)$, let $\bw^L(t)= (w_0^L, w_1^L,...,w^L_n, \fw_1^L,\ldots,\fw_n^L) \in\bV^L$ be the solution of the problem
\beq
B(\bw^L(t)-\bu(t), \bv^L) +A(\bw^L(t)-\bu(t),\bv^L)=0
\label{eq:eqforbwL}
\eeq
for all $\bv^L\in \bV^L$.
As the coefficients $a$ and $b$ in \eqref{eq:coercive} are both uniformly bounded and coercive for all $x\in D$ and $\by\in\bY$, problem \eqref{eq:eqforbwL} has a unique solution. 
Let $\boldq^L=\bw^L-\bu$. We then have the following estimate. 
\begin{lemma}
\label{lem:estqL}
For the solution  $\bw^L$ of problem \eqref{eq:eqforbwL}
\[
%\beq
\left\|\boldq^L(t)\right\|_{\bV} \leq c \inf_{\bv^L \in \bV^L} \left\|\bu(t)-\bv^L\right\|_{\bV}.
%\eeq
\]
\end{lemma}
\bproof 
From \eqref{eq:eqforbwL}, we have
\[
B(\bw^L-\bu,\bw^L-\bu)+A(\bw^L-\bu,\bw^L-\bu)=B(\bw^L-\bu,\bv^L-\bu)+A(\bw^L-\bu,\bv^L-\bu)
\]
for all $\bv^L \in \bV^L$. From the coerciveness and boundedness of the matrices $a$  and $b$ we get the conclusion. \eproof

When $\bu$ is sufficiently regular with respect to $t$, we have the following estimates.  

\begin{lemma}\label{lem:estdtqL}
If $\frac{\partial \bu}{\partial t}  \in  C([0,T], \bV)$, then
$$\left\|\frac{\partial \boldq^L}{\partial t}\right\|_{L^\infty(0,T;\bV)} \leq c \sup_{t \in [0,T]} \inf _{\bv^L\in \bV^L}\left\|\frac{\partial \bu}{\partial t}-\bv^L\right\|_{\bV}.$$
If $\frac{\partial^2 \bu}{\partial t^2}  \in  L^2(0,T, \bV)$, then
$$\left\|\frac{\partial^2 \boldq^L}{\partial t^2}\right\|_{L^2(0,T;\bV)} \leq c\inf _{\bv^L\in L^2(0,T;\bV^L)}\left\|\frac{\partial^2 \bu}{\partial t^2}-\bv^L\right\|_{L^2(0,T;\bV)}.$$
\end{lemma}

\bproof
If $\frac{\partial \bu}{\partial t}  \in  C([0,T]; \bV)$ from \eqref{eq:eqforbwL} we have
\[
B\left(\frac{\partial}{\partial t}\bw^L(t)-{\partial\over\partial t}\bu(t), \bv^L\right)+A\left(\frac{\partial}{\partial t}\bw^L(t)-{\partial\over\partial t}\bu(t), \bv^L\right) =0
\]
for all $\bv^L \in \bV^L$. We then proceed as in the proof of Lemma \ref{lem:estqL} to show the first inequality. The proof for the second inequality is similar.
\eproof

Let $\boldp^L=\bu^L-\bw^L$, i.e., for $i=1,\ldots,n$, 
$
p^L_i=u^L_i-w_i^L,  \frak{p}^L_i=\fu^L_i-\fw_i^L
 \text{ and }
p^L_0=u^L_0-w_0^L.
$
%We denote by $\boldp^L=(p_0^L,p_1^L,\ldots,p_n^L, \fp_1^L,\ldots,\fp_n^L)$. 
%
We recall the definition of the spaces $H_i$ in \eqref{eq:H}.
\bp
\label{prop:semidiscpest}
Assume that ${\partial^2\bu\over\partial t^2}\in L^2(0,T;\bV)$. Then there is a constant $c$ depending on $T$ such that for all $t\in (0,T)$
\begin{align*}
&\left\|{\partial p_0^L\over\partial t}(t)+\sum_{i=1}^n\nabla_{y_i}{\partial \fp_i^L\over\partial t}(t)\right\|_{H_n}+\left\|\curl p_0^L(t)+\sum_{i=1}^n\scurl_{y_i}p_i^L(t)\right\|_{H_n}\\
&\quad\le c\left[\left\|{\partial^2q_0^L\over\partial t^2}+\sum_{i=1}^n\nabla_{y_i}{\partial^2\fq_i^L\over\partial t^2}-q_0^L-\sum_{i=1}^n\nabla_{y_i}\fq^L_i\right\|_{L^2(0,T;H_n)}\right.\\
&\qquad~\left.+\left\|{\partial p_0^L\over\partial t}(0)+\sum_{i=1}^n\nabla_{y_i}{\partial \fp_i^L\over\partial t}(0)\right\|_{H_n}+\left\|\curl p_0^L(0)\right\|_H\right]. 
\end{align*}
\epr
\bproof
Since ${\partial^2\bu\over\partial t^2}\in L^2(0,T;\bV)$, from \eqref{eq:msprob} and \eqref{eq:nmsprob} we have for all $\bv^L=(v_0^L,v_1^L,\ldots,v_n^L,\fv_1^L,\ldots,\fv_n^L)\in\bV^L$
\begin{align*}
&\int_D\int_\bY \left[b(x,\by)\left({\partial^2 p_0^L\over\partial t^2}+\sum_{i=1}^n\nabla_{y_i}{\partial^2\fp_i^L\over\partial t^2}\right)\cdot\left(v_0^L+\sum_{i=1}^n\nabla_{y_i}\fv_i^L\right)\right.\\
&\qquad+\left.a(x,\by)\left(\curl p_0^L+\sum_{i=1}^n\scurl_{y_i}p_i^L\right)\cdot\left(\curl v_0^L+\sum_{i=1}^n\scurl_{y_i}v_i^L\right)\right]d\by dx\\
&\quad=-\!\int_D\int_\bY\!\! b(x,\by)\left({\partial^2q_0^L\over\partial t^2}+\sum_{i=1}^n\nabla_{y_i}{\partial^2\fq_i^L\over\partial t^2}\right)\!\cdot\!\left(v_0^L+\sum_{i=1}^n\nabla_{y_i}\fv_i^L\right)d\by dx-A(\boldq^L, \bv^L).
\end{align*}
From \eqref{eq:eqforbwL} we have $A(\boldq^L,\bv^L)=-B(\boldq^L,\bv^L)$. Thus
\begin{align}
&\int_D\int_\bY \left[b(x,\by)\left({\partial^2 p_0^L\over\partial t^2}+\sum_{i=1}^n\nabla_{y_i}{\partial^2\fp_i^L\over\partial t^2}\right)\cdot\left(v_0^L+\sum_{i=1}^n\nabla_{y_i}\fv_i^L\right)\right.\nonumber\\
&\qquad+\left.a(x,\by)\left(\curl p_0^L+\sum_{i=1}^n\scurl_{y_i}p_i^L\right)\cdot\left(\curl v_0^L+\sum_{i=1}^n\scurl_{y_i}v_i^L\right)\right]d\by dx\nonumber\\
&\quad=-\int_D\int_\bY b(x,\by)\left({\partial^2q_0^L\over\partial t^2}+\sum_{i=1}^n\nabla_{y_i}{\partial^2\fq_i^L\over\partial t^2}-q_0^L-\sum_{i=1}^n\nabla_{y_i}\fq^L_i\right)\nonumber\\
&\qquad\qquad\qquad\qquad\qquad\qquad\qquad\qquad\cdot\left(v_0^L+\sum_{i=1}^n\nabla_{y_i}\fv_i^L\right)d\by dx.
\label{eq:eqforpL}
\end{align}
Let $\bv^L={\partial\boldp^L\over\partial t}$. We then have 
\begin{align*}
&\frac12{d\over dt}\int_D\int_\bY \left[b(x,\by)\left({\partial p_0^L\over\partial t}+\sum_{i=1}^n\nabla_{y_i}{\partial \fp_i^L\over\partial t}\right)\cdot\left({\partial p_0^L\over\partial t}+\sum_{i=1}^n\nabla_{y_i}{\partial\fp_i^L\over\partial t}\right)\right.\\
&\qquad\left.+a(x,\by)\left(\curl p_0^L+\sum_{i=1}^n\scurl_{y_i}p_i^L\right)\cdot\left(\curl p_0^L+\sum_{i=1}^n\scurl_{y_i}p_i^L\right)\right]d\by dx\\
&\quad\le c~\left\|{\partial^2q_0^L\over\partial t^2}+\sum_{i=1}^n\nabla_{y_i}{\partial^2\fq^L_i\over\partial t^2}-q_0^L-\sum_{i=1}^n\nabla_{y_i}\fq^L_i\right\|_{H_n}\left\|{\partial p_0^L\over\partial t}+\sum_{i=1}^n\nabla_{y_i}{\partial \fp_i^L\over\partial t}\right\|_{H_n}\\
&\quad\le {c\over\gamma}\left\|{\partial^2q_0^L\over\partial t^2}+\sum_{i=1}^n\nabla_{y_i}{\partial^2\fq^L_i\over\partial t^2}-q_0^L-\sum_{i=1}^n\nabla_{y_i}\fq^L_i\right\|_{H_n}^2+c\gamma\left\|{\partial p_0^L\over\partial t}+\sum_{i=1}^n\nabla_{y_i}{\partial \fp_i^L\over\partial t}\right\|^2_{H_n}
\end{align*}
for a constant $\gamma>0$.
Integrating both sides on $(0,t)$ for $0<t<T$, and using the coercivity of the matrices $a$ and $b$, we have
\begin{align*}
&\left\|{\partial p_0^L\over\partial t}(t)+\sum_{i=1}^n\nabla_{y_i}{\partial \fp_i^L\over\partial t}(t)\right\|_{H_n}^2+\left\|\curl p_0^L(t)+\sum_{i=1}^n\scurl_{y_i}p_i^L(t)\right\|_{H_n}^2\\
&\quad\le 
{c\over\gamma}\left\|{\partial^2q_0^L\over\partial t^2}+\sum_{i=1}^n\nabla_{y_i}{\partial^2\fq^L_i\over\partial t^2}-q_0^L-\sum_{i=1}^n\nabla_{y_i}\fq^L_i\right\|^2_{L^2(0,T;H_n)}\\
&\qquad+c\gamma T\sup_{t\in [0,T]}\left\|{\partial p_0^L\over\partial t}(t)+\sum_{i=1}^n\nabla_{y_i}{\partial \fp_i^L\over\partial t}(t)\right\|^2_{H_n}\\
&\qquad+c\left\|{\partial p_0^L\over\partial t}(0)
+\sum_{i=1}^n\nabla_{y_i}{\partial \fp_i^L\over\partial t}(0)\right\|^2_{H_n}+c\left\|\curl p_0^L(0)+\sum_{i=1}^n\scurl_{y_i}p_i^L(0)\right\|^2_{H_n}.
\end{align*}
Choosing a sufficiently small constant $\gamma$, there is a constant $c$ depending on $T$ so that for all $t\in (0,T)$
\begin{align*}
&\left\|{\partial p_0^L\over\partial t}(t)+\sum_{i=1}^n\nabla_{y_i}{\partial \fp_i^L\over\partial t}(t)\right\|_{H_n}^2+\left\|\curl p_0^L(t)+\sum_{i=1}^n\scurl_{y_i}p_i^L(t)\right\|_{H_n}^2\\
&\quad\le c\left[\left\|{\partial^2q_0^L\over\partial t^2}+\sum_{i=1}^n\nabla_{y_i}{\partial^2\fq^L_i\over\partial t^2}-q_0^L-\sum_{i=1}^n\nabla_{y_i}\fq^L_i\right\|_{L^2(0,T;H_n)}^2\right.\\
&\qquad\left.+\left\|{\partial p_0^L\over\partial t}(0)+\sum_{i=1}^n\nabla_{y_i}{\partial \fp_i^L\over\partial t}(0)\right\|_{H_n}^2+\left\|\curl p_0^L(0)+\sum_{i=1}^n\scurl_{y_i}p_i^L(0)\right\|_{H_n}^2\right]. 
\end{align*}
Consider equation \eqref{eq:eqforpL} for $t=0$. Let $v_0^L=0$, $\fv_i^L=0$ and $v_i^L=p_i^L$. We then have 
\[
\int_D\int_\bY a(x,\by)\left(\curl p_0^L(0)+\sum_{i=1}^n\scurl_{y_i}p_i^L(0)\right)\cdot\left(\sum_{i=1}^n\scurl_{y_i}p_i^L(0)\right)d\by dx=0,
\]
i.e.,
\begin{multline*}
\int_D\int_\bY a(x,\by)\left(\sum_{i=1}^n\scurl_{y_i}p_i^L(0)\right)\cdot\left(\sum_{i=1}^n\scurl_{y_i}p_i^L(0)\right)d\by dx\\=-\int_D\int_\bY a(x,\by)\curl p_0^L(0)\cdot\left(\sum_{i=1}^n\scurl_{y_i}p_i^L(0)\right)d\by dx.
\end{multline*}
Using \eqref{eq:coercive}, we deduce that
\[
\left\|\sum_{i=1}^n\scurl_{y_i}p_i^L(0)\right\|_{H_n}\le c\left\|\curl p_0^L(0)\right\|_{H}.
\]
We then get the conclusion.\eproof
\begin{proposition}\label{prop:semidiscerror}

Assume that $\frac{\partial ^2 \bu }{\partial t^2} \in  L^2(0,T; \bV)$, and that
\be
\lim\limits_{L \to \infty} \|g_0^L-g_0\|_{W}=0 \ \ \mbox{and}\ \  \lim\limits_{L \to \infty} \|g_1^L-g_1\|_{H}=0.
\label{eq:inicondL}
\ee
 Then
\beqas
&&\lim\limits_{L \to \infty} \left \{  \left\|\frac{\partial (u_0^L-u_0)}{\partial t}\right\|_{L^\infty (0,T;H)}+ \sum_{i=1}^n\left\|\nabla_{y_i}\frac{ \partial (\frak u_i^L-\frak u_i)}{\partial t}\right\|_{L^\infty (0,T;H_i)}\right.\\
&&\qquad\qquad\left.+\left\|\curl (u_0^L-u_0)\right\|_{L^\infty (0,T;H)}+\sum_{i=1}^n\left\|\scurl_{y_i}(u_i^L-u_i)\right\|_{L^\infty(0,T;H_i)}\right \}=0.
\eeqas
\end{proposition}

\bproof
From Proposition \ref{prop:semidiscpest}, as $\bu^L-\bu=\boldp^L+\boldq^L$, we have
\begin{align}
&\left\|\frac{\partial (u_0^L-u_0)}{\partial t}+ \sum_{i=1}^n\nabla_{y_i}\frac{ \partial (\frak u_i^L-\frak u_i)}{\partial t}\right\|^2_{L^\infty (0,T; H_n)}\nonumber\\
&\qquad\qquad+\left\|\curl(u_0^L-u_0)+\sum_{i=1}^n\scurl_{y_i}(u^L_i-u_i)\right\|^2_{L^\infty (0,T;H_n)}\nonumber\\
&\qquad \leq c \left[\left\|{\partial^2q_0^L\over\partial t^2}+\sum_{i=1}^n\nabla_{y_i}{\partial^2\fq^L_i\over\partial t^2}-q_0^L-\sum_{i=1}^n\nabla_{y_i}\fq^L_i\right\|_{L^2(0,T;H_n)}^2\right.\nonumber\\
&\qquad\qquad\left.+\left\|{\partial p_0^L\over\partial t}(0)+\sum_{i=1}^n\nabla_{y_i}{\partial \fp_i^L\over\partial t}(0)\right\|_{H_n}^2+\left\|\curl p_0^L(0)\right\|_{H}^2\right]\nonumber\\
&\qquad\qquad +\left\|\frac{\partial q_0^L}{\partial t}+ \sum_{i=1}^n\nabla_{y_i}\frac{ \partial \fq_i^L}{\partial t}\right\|^2_{L^\infty (0,T; H_n)}+\|\boldq^L\|^2_{L^\infty (0,T;\bV)}
\label{eq:semidiscerrorest}
\end{align}
We show that $\lim\limits_{L \to \infty}\|\boldq^L\|_{L^\infty (0,T;\bV)}=0.$ As $\bu \in  C([0, T]; \bV), \bu  $ is uniformly continuous
as a function from $[0, T]$ to $\bV$. For $\varepsilon > 0$, there is a piecewise constant (with respect to $t$) function
 $\tilde{\bu} \in L^\infty (0,T; \bV)$ such that $\|\bu -\tilde{\bu} \|_{L^\infty (0,T; \bV)} < \varepsilon$. As $\tilde{\bu}(t)$ obtains only a finite number of $\bV$-values, when $L$ is sufficiently large, there is $\bv^L \in L^\infty (0,T; \bV^L)$ such that $\|\tilde{\bu}-\bv^L \|_{L^\infty (0,T; \bV)} < \varepsilon$. Thus
 $$\lim_{L \to \infty} \sup_{t\in (0,T)} \inf _{\bv^L\in \bV^L} \|\bu(t)-\bv^L\|_{\bV}=0. $$
 We then apply Lemma \ref{lem:estqL}.
 Similarly, we have from Lemmas \ref{lem:estqL} and \ref{lem:estdtqL}
\[
\lim_{L\to\infty}\left\|{\partial^2q_0^L\over\partial t^2}+\sum_{i=1}^n\nabla_{y_i}{\partial^2\fq^L_i\over\partial t^2}-q_0^L-\sum_{i=1}^n\nabla_{y_i}\fq^L_i\right\|_{L^2(0,T;H_n)}=0
\]
and 
\[
\lim_{L \to \infty}\left\|\frac{\partial q_0^L}{\partial t}+ \sum_{i=1}^n\nabla_{y_i}\frac{ \partial \fq^L_i}{\partial t}\right\|_{L^\infty (0,T; H_n)}=0.
\]
Furthermore, we have that
$$\left\|\curl p_0^L(0)\right\|_{H} \le \left\|\curl u_0^L(0)-\curl u_0(0)\right\|_{H}+\left\|\curl u_0(0)-\curl w_0^L(0)\right\|_{H},$$
which converges to $0$ due to \eqref{eq:inicondL} and Lemma \ref{lem:estqL}. Similarly, we have
$$\lim_{L \to \infty}\left\|\frac{\partial  p^L_0}{\partial t}(0)+\sum_{i}\nabla_{y_i}\frac{\partial \fp_i^L}{\partial t}(0)\right\|_{H_n}=0.$$
We then get the conclusion.\eproof
%%%%%%%%%%%%%%%%%%%%%%%%%%%%%%%%%%%%%%%%%%%%%%%%%%%%%%%%%%%%%%%%%%%%%%
\subsection{Fully discrete problem}
Following the scheme of Dupont \cite{Dupont}, we discretize problem \eqref{eq:nmsprob} in both spatial and temporal variables.
Let $\Delta t= \frac{T}{M}$ where $M$ is a positive integer. Let $t_m= m \Delta t$. We employ the following notations of Dupont for a function $r \in C([0,T]; X)$ where $X$ is a Banach space and $r_m= r(t_m, \cdot)$
\begin{align*}
    r_{m+1/2}&=\frac{1}{2}(r_{m+1}+r_m),& r_{m,\theta}= \theta r_{m+1}+(1-2\theta)r_m+\theta r_{m-1},\\
    \partial _t  r_{m+1/2}&=(r_{m+1}-r_m)/\Delta t,& \partial^2_t r_m=(r_{m+1}-2 r_m+r_{m-1})/ (\Delta t)^2,\\
    \delta _t r_m&= (r_{m+1}-r_{m-1})/ (2 \Delta t).&
\end{align*}

We consider the following fully discrete problem: 

For $m=1,...,M$ find $\bu^L_m=(u_{0,m}^L,u_{1,m}^L,...,u_{n,m}^L, \frak u_{1,m}^L ,...,\frak u_{n,m}^L) \in \bV^L$ such that for $m=1,...,M-1$
\begin{eqnarray}
&&\int_D\int_\bY\left[b(x,\by)\left(\partial^2_tu^L_{0,m}+\sum_{i=1}^n\nabla_{y_i} \partial^2_t\frak u^L_{i,m}\right)\cdot\left(v_0^L+\sum_{i=1}^n\nabla_{y_i} \frak v^L_i\right)+\right.\nonumber\\
&&\left.a(x,\by)\left(\curl u^L_{0 ,m , 1/4}+\sum_{i=1}^n\scurl_{y_i}u^L_{i ,m , 1/4}\right)\cdot\left(\curl v^L_0+\sum_{i=1}^n\scurl_{y_i}v^L_i\right)\right]d\by dx \nonumber\\
&&=\int_D f_{m,1/4}(t,x)\cdot v^L_0(x)dx,
\label{eq:fullydiscprob}
\end{eqnarray}
for all $\bv^L=(v_0^L,v_1^L,\ldots,v_n^L,\fv_1^L,\ldots,\fv_n^L) \in \bV^L$. 

For continuous functions $r: [0,T] \to X$, let
$$\|r\|_{\tilde{L}^\infty (0,T;X)}:= \max_{0 \leq m < M}\|r_{m+1/2}\|_{X}.$$
We also denote by
$$\|\partial _t r\|_{\tilde{L}^\infty (0,T;X)}:= \max_{0 \leq m < M}\|\partial_tr_{m+1/2}\|_{X}.$$
Let
$$\boldp^L_m:= \bu^L_m- \bw^L_m.$$
\begin{lemma}
\label{lem:fullydiscpest}
Assume that $\bu \in H^2(0,T; \bV)$, ${\partial^2q^L_0\over\partial t^2}\in L^2(0,T;H)$, $\frac{\partial ^2}{\partial t^2}\nabla_{y_i}\fq^L_{i} \in L^2(0,T;H_i)$. If $\frac{\partial^3u_0}{\partial t^3}\in L^2(0,T; H)$ and $\frac{\partial^3}{\partial t^3}\nabla_{y_i} \frak u_i \in L^2(0,T; H_i)$, then there exists a constant $c$ independent of $\Delta t$ and $\bu$
such that for each $j=1,2,...,M-1$
\begin{align*}
&\|\partial_t p^L_{0,j+1/2}\|^2_H+\sum_{i=1}^n\|\partial_t\nabla_{y_i}\fp_{i,j+1/2}^L\|^2_{H_i}+\|\curl p^L_{0,j+1/2}\|_H^2+\sum_{i=1}^n\|\scurl_{y_i}p^L_{i,j+1/2}\|_{H_i}^2\\
&\quad\le c\left[(\Delta t)^2 \left\|{\partial^3u_0\over\partial t^3}\right\|_{H}^2+(\Delta t)^2\sum_{i=1}^n\left\|{\partial^3\nabla_{y_i}\fu_i\over\partial t^3}\right\|^2_{H_i}+\left\|{\partial^2q^L_0\over\partial t^2}\right\|^2_{L^2(0,T;H)}\right.\\
&\qquad\left.+\sum_{i=1}^n\left\|{\partial^2\over\partial t^2}\nabla_{y_i}\fq_{i}^L\right\|_{L^2(0,T;H_i)}^2+\|q_0^L\|_{L^\infty(0,T;H)}^2+\sum_{i=1}^n\|\nabla_{y_i}\fq^L_i\|_{L^\infty(0,T;H_i)}^2\right]\\
&\qquad+c\left(\|\partial_t p^L_{0,1/2}\|^2_{H}+\sum_{i=1}^n\|\partial_t\nabla_{y_i} \fp^L_{i,1/2}\|_{H_i}^2+\|\curl p^L_{0,1/2}\|^2_H+\sum_{i=1}^n\|\scurl_{y_i}p^L_{i,1/2}\|^2_{H_i}\right).
\end{align*}
Further, if $\frac{\partial^4u_0}{\partial t^4}\in L^2(0,T;H)$ and $\frac{\partial^4}{\partial t^4}\nabla_{y_i} \frak u_i \in L^2(0,T; H_i)$, then there exists a constant $c$ independent of $\Delta t$ and $\bu$ such that for each $j=1,2,...,M-1$
\begin{align*}
&\|\partial_t p^L_{0,j+1/2}\|^2_H+\sum_{i=1}^n\|\partial_t\nabla_{y_i}\fp_{i,j+1/2}^L\|^2_{H_i}+\|\curl p^L_{0,j+1/2}\|^2_H+\sum_{i=1}^n\|\scurl_{y_i}p^L_{i,j+1/2}\|^2_{H_i}\\
&\quad\le c\left[(\Delta t)^4 \left\|{\partial^4u_0\over\partial t^4}\right\|_{H}^2+(\Delta t)^4\sum_{i=1}^n\left\|{\partial^4\nabla_{y_i}\fu_i\over\partial t^4}\right\|^2_{H_i}+\left\|{\partial^2q^L_0\over\partial t^2}\right\|^2_{L^2(0,T;H)}\right.\\
&\qquad\left.+\sum_{i=1}^n\left\|{\partial^2\over\partial t^2}\nabla_{y_i}\fq_{i}^L\right\|_{L^2(0,T;H_i)}^2+\|q_0^L\|_{L^\infty(0,T;H)}^2+\sum_{i=1}^n\|\nabla_{y_i}\fq^L_i\|_{L^\infty(0,T;H_i)}^2\right]\\
&\qquad+c\left(\|\partial_t p^L_{0,1/2}\|^2_{H}+\sum_{i=1}^n\|\partial_t\nabla_{y_i} \fp^L_{i,1/2}\|_{H_i}^2+\|\curl p^L_{0,1/2}\|^2_H+\sum_{i=1}^n\|\scurl_{y_i}p^L_{i,1/2}\|^2_{H_i}\right).
\end{align*}
\end{lemma}
\bproof
From \eqref{eq:eqforbwL} and \eqref{eq:fullydiscprob}, we have
\beqas
&&A(\bw^L,\bv^L)=A(\bu,\bv^L)-B(\bw^L-\bu,\bv^L)=\int_D f(t,x)\cdot v_0^L(x)dx\\
&&-\int_D\int_\bY b(x,\by)\left({\partial^2u_0\over\partial t^2}+\sum_{i=1}^n{\partial^2\over\partial t^2}\nabla_{y_i}\fu_i\right)\cdot\left(v_0^L+\sum_{i=1}^n\nabla_{y_i}\fv_i^L\right)d\by dx -B(\boldq^L, \bv^L).
\eeqas
Averaging this equation at  $t_{m+1}$, $t_m$ and $t_{m-1}$ with weights $\frac{1}{4}, \frac{1}{2}, \frac{1}{4}$ respectively, and using \eqref{eq:fullydiscprob}, we get
\begin{align*}
&\int_D\int_\bY b(x,\by)\left(\partial^2_tu_{0,m}^L+\sum_{i=1}^n\nabla_{y_i}\partial^2_t\fu_{i,m}^L\right)\cdot\left(v_0^L+\sum_{i=1}^n\nabla_{y_i}\fv_i^L\right)d\by dx+A(\boldp^L_{m,1/4},\bv^L)\\
&\quad=\int_D\int_\bY b(x,\by)\left({\partial^2 u_{0,m,1/4}\over\partial t^2}+\sum_{i=1}^n{\partial^2\over\partial t^2}\nabla_{y_i}\fu_{i,m,1/4}\right)\cdot\left(v_0^L+\sum_{i=1}^n\nabla_{y_i}\fv_i^L\right)d\by dx\\
&\qquad+B(\boldq^L_{m,1/4},\bv^L).
\end{align*}
Thus
\begin{align*}
&\int_D\int_\bY b(x,\by)\left(\partial^2_tp^L_{0,m}+\sum_{i=1}^n\nabla_{y_i}\partial^2_t\fp_{i,m}^L\right)\cdot\left(v_0^L+\sum_{i=1}^n\nabla_{y_i}\fv_i^L\right)d\by dx+A(\boldp_{m,1/4}^L,\bv^L)\\
&\quad=\int_D\int_\bY b(x,\by)\left({\partial^2 u_{0,m,1/4}\over\partial t^2}-\partial^2_tu_{0,m}+\sum_{i=1}^n\left({\partial^2\over\partial t^2}\nabla_{y_i}\fu_{i,m,1/4}-\nabla_{y_i}\partial^2_t\fu_{i,m}\right)\right)\\
&\qquad\qquad\qquad\qquad\qquad\qquad\qquad\qquad\qquad\qquad\qquad\qquad\cdot\left(v_0^L+\sum_{i=1}^n\nabla_{y_i}\fv_i^L\right)d\by dx\\
&\qquad-\int_D\int_\bY b(x,\by)\left(\partial^2_t q_{0,m}^L+\sum_{i=1}^n\nabla_{y_i}\partial^2_t q_{i,m}^L\right)\cdot\left(v_0^L+\sum_{i=1}^n\nabla_{y_i}\fv_i^L\right)d\by dx\\
&\qquad+B(\boldq^L_{m,1/4},\bv^L).
\end{align*}
We denote by
\[
s_{0,m}={\partial^2 u_{0,m,1/4}\over\partial t^2}-\partial^2_t u_{0,m},\ \ s_{i,m}={\partial^2\over\partial t^2}\nabla_{y_i}\fu_{i,m,1/4}-\partial^2_t\nabla_{y_i}\fu_{i,m}.
\]
Let $\bv^L=\delta_t\boldp^L_m$. Using the following relationships:
\begin{align*}
\partial^2_t r_m=& \frac{1}{\Delta t} (\partial_t r_{m+1/2}-\partial _t r_{m-1/2}), \qquad\qquad r_{m, 1/4}= \frac{1}{2}(r_{m+1/2}+r_{m-1/2})\\
\delta_t r_m=& \frac{1}{2}(\partial _t r_{m+1/2}+\partial _t r_{m-1/2})=\frac{1}{\Delta t}(r_{m+1/2}-r_{m-1/2}),
\end{align*}
we have
\begin{align*}
&{1\over 2\Delta t}\int_D\int_\bY b(x,\by)\left(\partial_t p^L_{0,m+1/2}-\partial_tp^L_{0,m-1/2}+\sum_{i=1}^n\nabla_{y_i}\left(\partial_t\fp_{i,m+1/2}^L-\partial_t\fp_{i,m-1/2}^L\right)\right)\\
&\qquad\qquad\cdot\left(\partial_t p^L_{0,m+1/2}+\partial_tp^L_{0,m-1/2}+\sum_{i=1}^n\nabla_{y_i}\left(\partial_t\fp_{i,m+1/2}^L+\partial_t\fp_{i,m-1/2}^L\right)\right)d\by dx
\end{align*}
\begin{align*}
&+{1\over 2\Delta t}\!\!\int_D\!\int_\bY \!\!a(x,\by)\left(\curl\left(p^L_{0,m+1/2}+p^L_{0,m-1/2}\right)+\!\sum_{i=1}^n\scurl_{y_i}\left(p^L_{i,m+1/2}+p^L_{i,m-1/2}\right)\right)\\
&\qquad\qquad\cdot\left(\curl\left(p^L_{0,m+1/2}-p^L_{0,m-1/2}\right)+\sum_{i=1}^n\scurl_{y_i}\left(p^L_{i,m+1/2}-p^L_{i,m-1/2}\right)\right)d\by dx\\
&=\frac12\!\int_D\!\int_\bY \!\!b(x,\by)\left(s_{0,m}-\partial^2_tq_{0,m}^L+q_{0,m,1/4}^L+\!\!\sum_{i=1}^n\left(s_{i,m}-\nabla_{y_i}\partial^2_t\fq_{i,m}^L+\nabla_{y_i}\fq_{i,m,1/4}^L\right)\right)\\
&\qquad\qquad\cdot\left(\partial_t p^L_{0,m+1/2}+\partial_t p^L_{0,m-1/2}+\sum_{i=1}^n\left(\nabla_{y_i}\partial_t \fp^L_{i,m+1/2}+\nabla_{y_i}\partial_t \fp^L_{i,m-1/2}\right)\right)d\by dx.
\end{align*}
We thus have
\begin{align*}
&{1\over 2\Delta t}\left[B\left(\partial_t\boldp^L_{m+1/2},\partial_t\boldp^L_{m+1/2}\right)-B\left(\partial_t\boldp^L_{m-1/2},\partial_t\boldp^L_{m-1/2}\right)\right.\\
&\qquad\left.+A\left(\boldp^L_{m+1/2},\boldp^L_{m+1/2}\right)-
A\left(\boldp^L_{m-1/2},\boldp^L_{m-1/2}\right)\right]\\
&\quad\le c\|s_{0,m}-\partial^2_tq^L_{0,m}+q^L_{0,m,1/4}+\sum_{i=1}^n(s_{i,m}-\nabla_{y_i}\partial^2_t\fq^L_{i,m}+\nabla_{y_i}\fq^L_{i,m,1/4})\|_{H_n}\\
&\qquad\cdot\left\|\partial_tp^L_{0,m+1/2}+\partial_tp^L_{0,m-1/2}+\sum_{i=1}^n\left(\nabla_{y_i}\partial_t\fp^L_{i,m+1/2}+\nabla_{y_i}\partial_t\fp^L_{i,m-1/2}\right)\right\|_{H_n}\\
&\quad\le {c\over\gamma}\left(\|s_{0,m}\|_{H}^2+\sum_{i=1}^n\|s_{i,m}\|_{H_i}^2+\|\partial^2_tq^L_{0,m}\|_H^2+\sum_{i=1}^n\|\nabla_{y_i}\partial^2_t\fq_{i,m}^L\|_{H_i}^2\right.\\
&\qquad\left.+\|q_{0,m,1/4}^L\|_H^2+\sum_{i=1}^n\|\nabla_{y_i}\fq^L_{i,m,1/4}\|^2_{H_i}\right)\\
&\qquad+ c\gamma\Big(\|\partial_tp^L_{0,m+1/2}\|_H^2+\|\partial_tp^L_{0,m-1/2}\|_H^2\\
&\qquad+\sum_{i=1}^n\|\nabla_{y_i}\partial_t\fp^L_{i,m+1/2}\|_{H_i}^2+\sum_{i=1}^n\|\nabla_{y_i}\partial_t\fp^L_{i,m-1/2}\|_{H_i}^2\Big).
\end{align*}
Summing this up for all $m=1, \ldots,j$, we deduce
\begin{align*}
&B(\partial_t\boldp^L_{j+1/2},\partial_t\boldp^L_{j+1/2})-B(\partial_t\boldp^L_{1/2},\partial_t\boldp^L_{1/2})+A(\boldp^L_{j+1/2},\boldp^L_{j+1/2})-A(\boldp^L_{1/2},\boldp^L_{1/2})\\
&\quad\le {c\over\gamma}{2\Delta t}\sum_{m=1}^M\left(\|s_{0,m}\|_H^2+\sum_{i=1}^n\|s_{i,m}\|_{H_i}^2+\|\partial^2_tq^L_{0,m}\|_H^2+\sum_{i=1}^n\|\nabla_{y_i}\partial^2_t\fq_{i,m}^L\|_{H_i}^2\right.\\
&\qquad\qquad\qquad\qquad\left.+\|q_{0,m,1/4}^L\|_H^2+\sum_{i=1}^n\|\nabla_{y_i}\fq^L_{i,m,1/4}\|_{H_i}^2\right)\\
&\qquad+c\gamma4\Delta tM\left(\max_{1\le m\le M}\|\partial_tp^L_{0,m+1/2}\|^2_{H}+\sum_{i=1}^n\max_{1\le m\le M}\|\partial_t\nabla_{y_i}\fp_{i,m+1/2}^L\|^2_{H_i}\right)\\
&\qquad+c\gamma2\Delta t\left(\|\partial_t p^L_{0,1/2}\|^2_{H}+\sum_{i=1}^n\|\partial_t\nabla_{y_i} \fp^L_{i,1/2}\|_{H_i}^2\right).
\end{align*}
From \eqref{eq:coercive}, we have
\begin{align*}
&\|\partial_t p^L_{0,j+1/2}\|_H^2+\sum_{i=1}^n\|\partial_t\nabla_{y_i}\fp_{i,j+1/2}^L\|_{H_i}^2+\|\curl p^L_{0,j+1/2}\|_H^2+\sum_{i=1}^n\|\scurl_{y_i}p^L_{i,j+1/2}\|_{H_i}^2\\
&\quad\le {c\over\gamma}{2\Delta t}\sum_{m=1}^M\left(\|s_{0,m}\|_H^2+\sum_{i=1}^n\|s_{i,m}\|_{H_i}^2+\|\partial^2_tq^L_{0,m}\|_H^2+\sum_{i=1}^n\|\nabla_{y_i}\partial^2_t\fq_{i,m}^L\|_{H_i}^2\right.\\
&\qquad\left.+\|q_{0,m,1/4}^L\|_H^2+\sum_{i=1}^n\|\nabla_{y_i}\fq^L_{i,m,1/4}\|_{H_i}^2\right)\\
&\qquad+c\gamma 4\Delta tM\left(\max_{1\le m\le M}\|\partial_tp^L_{0,m+1/2}\|_{H}^2+\sum_{i=1}^n\max_{1\le m\le M}\|\partial_t\nabla_{y_i}\fp_{i,m+1/2}^L\|_{H_i}^2\right)\\
&\qquad+c\left(\|\partial_t p^L_{0,1/2}\|^2_{H}+\sum_{i=1}^n\|\partial_t\nabla_{y_i} \fp^L_{i,1/2}\|_{H_i}^2+\|\curl p^L_{0,1/2}\|_H^2+\sum_{i=1}^n\|\scurl_{y_i}p^L_{i,1/2}\|_{H_i}^2\right).
\end{align*}
Choosing $\gamma$ sufficiently small, we deduce that
\beqas
\|\partial_t p^L_{0,j+1/2}\|_H^2+\sum_{i=1}^n\|\partial_t\nabla_{y_i}\fp_{i,j+1/2}^L\|_{H_i}^2+\|\curl p^L_{0,j+1/2}\|_H^2+\sum_{i=1}^n\|\scurl_{y_i}p^L_{i,j+1/2}\|_{H_i}^2\\
\le {c\over\gamma}{2\Delta t}\sum_{m=1}^M\left(\|s_{0,m}\|_H^2+\sum_{i=1}^n\|s_{i,m}\|_{H_i}^2+\|\partial^2_tq^L_{0,m}\|_H^2+\sum_{i=1}^n\|\nabla_{y_i}\partial^2_t\fq_{i,m}^L\|_{H_i}^2\right.\\
\left.+\|q_{0,m,1/4}^L\|_H^2+\sum_{i=1}^n\|\nabla_{y_i}\fq^L_{i,m,1/4}\|_{H_i}^2\right)\\
+c\left(\|\partial_t p^L_{0,1/2}\|^2_{H}+\sum_{i=1}^n\|\partial_t\nabla_{y_i} \fp^L_{i,1/2}\|_{H_i}^2+\|\curl p^L_{0,1/2}\|_H^2+\sum_{i=1}^n\|\scurl_{y_i}p^L_{i,1/2}\|_{H_i}^2\right).
\eeqas
Following Dupont \cite{Dupont}, using the integral formula of the remainder of Taylor expansion, we have,
\begin{align*}
\partial^2_tq^L_{0,m}=
%{q_{0,m+1}^L-2q^L_{0,m}+q^L_{0,m-1}\over (\Delta t)^2}\\
%&={1\over(\Delta t)^{2}}\left(q_{0,m+1}^L-q^L_{0,m}-{\partial q_{0,m}^L\over \partial t}\Delta t+q_{0,m-1}^L-q^L_{0,m}+{\partial q_{0,m}^L\over \partial t}\Delta t\right)\\
%&={1\over(\Delta t)^{2}}\left(q_{0,m+1}^L-q^L_{0,m}-{\partial q_{0,m}^L\over \partial t}(t_{m+1}-t_{m})+q_{0,m-1}^L-q^L_{0,m}-{\partial q_{0,m}^L\over \partial t}(t_{m-1}-t_{m})\right)\\
%&={1\over(\Delta t)^{2}}\left(\int_{t_m}^{t_{m+1}} (t_{m+1}-t) {\partial^2q^L_{0}\over\partial t^2}(t)dt+\int_{t_m}^{t_{m-1}} (t_{m-1}-t) {\partial^2q^L_{0}\over\partial t^2}(t)dt\right)\\
%&={1\over(\Delta t)^{2}}\left(\int_{0}^{\Delta t} (\Delta t-\tau) {\partial^2q^L_{0}\over\partial t^2}(t_m+\tau)d\tau+\int_{0}^{-\Delta t} (-\Delta t-\tau) {\partial^2q^L_{0}\over\partial t^2}(t_m+\tau)d\tau\right)\\
(\Delta t)^{-2}\int_{-\Delta t}^{\Delta t}(\Delta t-|\tau|){\partial^2q^L_{0}\over\partial t^2}(t_m+\tau)d\tau,
\end{align*}
and similarly, for $i=1,\ldots,n$
\[
\partial ^2_t (\nabla_{y_i}\fq^L_{i, m})= (\Delta t)^{-2} \int_{-\Delta t}^{\Delta t} (\Delta t -|\tau|) \frac{\partial ^2 \nabla_{y_i}\fq^L_{i}}{\partial t^2}(t_m+\tau) d\tau.
\]
Using Cauchy-Schwarz inequality, we have
\begin{align*}
\sum_{m=1}^M\|\partial^2_tq^L_{0,m}\|^2_{H}\Delta t 
%&=\sum_{m=1}^M \int_D|\partial^2_tq^L_{0,m}|^2 dx\Delta t\\
%&=\sum_{m=1}^M \int_D\left|(\Delta t)^{-2}\int_{-\Delta t}^{\Delta t}(\Delta t-|\tau|){\partial^2q^L_{0}\over\partial t^2}(t_m+\tau)d\tau\right|^2 dx\Delta t\\
%&\le\sum_{m=1}^M (\Delta t)^{-3}\int_{-\Delta t}^{\Delta t}(\Delta t-|\tau|)^2d\tau\int_D\int_{-\Delta t}^{\Delta t}\left({\partial^2q^L_{0}\over\partial t^2}(t_m+\tau)\right)^2d\tau dx\\
%&\le\sum_{m=1}^M (\Delta t)^{-3} \frac23 (\Delta t)^{3}\int_D\int_{-\Delta t}^{\Delta t}\left({\partial^2q^L_{0}\over\partial t^2}(t_m+\tau)\right)^2d\tau dx\\
%&\le \frac23 \sum_{m=1}^M\int_D\int_{t_{m-1}}^{t_{m+1}}\left({\partial^2q^L_{0}\over\partial t^2}(t)\right)^2dt dx\\
\le \frac43\left\|{\partial^2q^L_{0}\over\partial t^2}\right\|^2_{L^2(0,T;H)},
\end{align*}
and similarly, we have
\[\sum_{m=1}^M\|\partial^2_t\nabla_{y_i}\fq^L_{i,m}\|^2_{H_i} \Delta t\le\frac43\left\|{\partial^2\over\partial t^2}\nabla_{y_i}\fq^L_{i}\right\|^2_{L^2(0,T;H_i)}.
\]
We write
\begin{align*}
 s_{0,m}=& \frac{1}{4}\int_0^{\Delta t} \left(1-2\left(1- \frac{|\tau|}{\Delta t}\right)^2\right) \frac{\partial ^3u_0}{\partial t^3}(t_m+\tau) d\tau-\frac{1}{4}\int^0_{-\Delta t} \left(1-2\left(1- \frac{|\tau|}{\Delta t}\right)^2\right) \frac{\partial^3u_0}{\partial t^3}(t_m+\tau) d\tau
\end{align*}
and
\begin{align*}
    s_{i,m}=& \frac{1}{4}\int_0^{\Delta t} \left(1-2\left(1- \frac{|\tau|}{\Delta t}\right)^2\right) \frac{\partial^3 \nabla_{y_i}\frak u_i}{\partial t^3}(t_m+\tau) d\tau
   -\frac{1}{4}\int^0_{-\Delta t} \left(1-2\left(1- \frac{|\tau|}{\Delta t}\right)^2\right) \frac{\partial^3\nabla_{y_i} \frak u_i}{\partial t^3}(t_m+\tau) d\tau.
\end{align*}
%Using similar estimate as above, we have
Therefore
\[
\|s_{0,m}\|^2_{H}\le c\Delta t\int_{t_{m-1}}^{t_m+1}\left\|{\partial^3u_0\over\partial t^3}(\tau)\right\|_H^2d\tau,\ \ \|s_{i,m}\|^2_{H_i}\le c \Delta t \int_{t_{m-1}}^{t_m+1}\left\|\frac{\partial ^3 \nabla_{y_i} \frak u_i}{\partial t^3}(\tau) \right\|_{H_i}^2 d\tau.
\]
We also have
\[
\|q^L_{0,m,1/4}\|_{H}\le \max_{t\in [0,T]}\|q^L_{0}(t)\|_{H}
\ \ 
\mbox{and}
\ \  
\|\nabla_{y_i}\fq^L_{i,m,1/4}\|_{H_i}\le \max_{t\in[0,T]} \|\nabla_{y_i}\fq^L_i(t)\|_{H_i}.
\]
We thus deduce
\begin{align*}
&\|\partial_t p^L_{0,j+1/2}\|_H^2+\sum_{i=1}^n\|\partial_t\nabla_{y_i}\fp_{i,j+1/2}^L\|^2_{H_i}+\|\curl p^L_{0,j+1/2}\|^2_H+\sum_{i=1}^n\|\scurl_{y_i}p^L_{i,j+1/2}\|^2_{H_i}\\
&\quad\le c\left[(\Delta t)^2 \left\|{\partial^3u_0\over\partial t^3}\right\|_H^2+(\Delta t)^2\sum_{i=1}^n\left\|{\partial^3\nabla_{y_i}\fu_i\over\partial t^3}\right\|^2_{H_i}+\left\|{\partial^2q^L_0\over\partial t^2}\right\|^2_{L^2(0,T;H)}\right.\\
&\qquad\left.+\sum_{i=1}^n\left\|{\partial^2\over\partial t^2}\nabla_{y_i}\fq_{i}^L\right\|_{L^2(0,T;H_i)}^2+\|q_0^L\|_{L^\infty(0,T;H)}^2+\sum_{i=1}^n\|\nabla_{y_i}\fq^L_i\|_{L^\infty(0,T;H_i)}^2\right]\\
&\qquad+c\left(\|\partial_t p^L_{0,1/2}\|^2_{H}+\sum_{i=1}^n\|\partial_t\nabla_{y_i} \fp^L_{i,1/2}\|_{H_i}^2+\|\curl p^L_{0,1/2}\|_H^2+\sum_{i=1}^n\|\scurl_{y_i}p^L_{i,1/2}\|^2_{H_i}\right).
\end{align*}
When
\[
{\partial^4u_0\over\partial t^4}\in L^2(0,T;H)\ \ \mbox{and}\ \ {\partial^4\over\partial t^4}\nabla_{y_i}\fu_i\in L^2(0,T;H_i),
\]
we have
\[
s_{0,m}=\frac{1}{12} \int_{- \Delta t}^{ \Delta t} (\Delta t - |\tau|)\left(3-2 \left(1-\frac{|\tau|}{\Delta t}\right)^2\right)\frac{\partial ^4u_0}{\partial t^4}(t_m+\tau) d \tau.
\]
and
\[  
s_{i,m}=\frac{1}{12} \int_{- \Delta t}^{ \Delta t} (\Delta t - |\tau|)\left(3-2 \left(1-\frac{|\tau|}{\Delta t}\right)^2\right)\frac{\partial ^4\nabla_{y_i} \frak u_i}{\partial t^4}(t_m+\tau) d \tau .
\]
%Using similar estimate as for $q_{0,m}^L$, we have
Therefore
\[
\|s_{0,m}\|^2_{H}\le c(\Delta t)^3\int_{t_{m-1}}^{t_{m+1}}\left\|{\partial^4u_0\over\partial t^4}(\tau)\right\|^2_{H}d\tau
\]
and
\[
\|s_{i,m}\|^2_{H_i}\leq c (\Delta t)^3 \int_{t_{m-1}}^{t_{m+1}}\left\|\frac{\partial ^4\nabla_{y_i} \frak u_i}{\partial t^4}(\tau) \right\|^2_{H_i}d\tau.
\]
Thus we have
\begin{align*}
&\|\partial_t p^L_{0,j+1/2}\|_H^2+\sum_{i=1}^n\|\partial_t\nabla_{y_i}\fp_{i,j+1/2}^L\|_{H_i}^2+\|\curl p^L_{0,j+1/2}\|_H^2+\sum_{i=1}^n\|\scurl_{y_i}p^L_{i,j+1/2}\|_{H_i}^2\\
&\quad\le c\left[(\Delta t)^4 \left\|{\partial^4u_0\over\partial t^4}\right\|_{H}^2+(\Delta t)^4\sum_{i=1}^n\left\|{\partial^4\nabla_{y_i}\fu_i\over\partial t^4}\right\|^2_{H_i}+\left\|{\partial^2q^L_0\over\partial t^2}\right\|^2_{L^2(0,T;H)}\right.\\
&\qquad\left.+\sum_{i=1}^n\left\|{\partial^2\over\partial t^2}\nabla_{y_i}\fq_{i}^L\right\|_{L^2(0,T;H_i)}^2+\|q_0^L\|_{L^\infty(0,T;H)}^2+\sum_{i=1}^n\|\nabla_{y_i}\fq^L_i\|_{L^\infty(0,T;H_i)}^2\right]\\
&\qquad+c\left(\|\partial_t p^L_{0,1/2}\|^2_{H}+\sum_{i=1}^n\|\partial_t\nabla_{y_i} \fp^L_{i,1/2}\|_{H_i}^2+\|\curl p^L_{0,1/2}\|_H^2+\sum_{i=1}^n\|\scurl_{y_i}p^L_{i,1/2}\|_{H_i}^2\right).
\end{align*}
\eproof

We then have the following error estimates.
\begin{proposition}
Assume that $\bu\in H^2(0,T;\bV)$. If ${\partial^3u_0\over\partial t^3}\in L^2(0,T;H)$ and ${\partial^3\nabla_{y_i}\fu_i\over\partial t^3}\in L^2(0,T;H_i)$, then there is a constant $c$ such that
\begin{align*}
&\|\partial_t u^L_0-\partial_tu_0\|_{\tL(0,T;H)}+\sum_{i=1}^n\|\partial_t\nabla_{y_i}\fu_i^L-\partial_t\nabla_{y_i}\fu_i\|_{\tL(0,T;H_i)}\\
&\qquad+\|\curl u^L_{0}-\curl u_0\|_{\tL(0,T;H)}+\sum_{i=1}^n\|\scurl_{y_i}u^L_{i}-\scurl_{y_i}u_i\|_{\tL(0,T;H_i)}\\
&\quad\le c\left[\Delta t \left\|{\partial^3u_0\over\partial t^3}\right\|_{L^2(0,T;H)}+\Delta t\sum_{i=1}^n\left\|{\partial^3\nabla_{y_i}\fu_i\over\partial t^3}\right\|_{L^2(0,T;H_i)}+\left\|{\partial^2q^L_0\over\partial t^2}\right\|_{L^2(0,T;H)}\right.\\
&\qquad\left.+\sum_{i=1}^n\left\|{\partial^2\over\partial t^2}\nabla_{y_i}\fq_{i}^L\right\|_{L^2(0,T;H_i)}+\|q_0^L\|_{L^\infty(0,T;H)}+\sum_{i=1}^n\|\nabla_{y_i}\fq^L_i\|_{L^\infty(0,T;H_i)}\right]\\
&\qquad+c\left(\|\partial_t p^L_{0,1/2}\|_{H}+\sum_{i=1}^n\|\partial_t\nabla_{y_i} \fp^L_{i,1/2}\|_{H_i}+\|\curl p^L_{0,1/2}\|_H+\sum_{i=1}^n\|\scurl_{y_i}p^L_{i,1/2}\|_{H_i}\right)\\
&\qquad+\|\partial_tq_0^L\|_{\tL(0,T;H)}+\|\curl q_0^L\|_{\tL(0,T;H)}\\
&\qquad+\sum_{i=1}^n\left(\|\partial_t\nabla_{y_i}\fq_i^L\|_{\tL(0,T;H_i)}+\|\scurl_{y_i}q_i^L\|_{\tL(0,T;H_i)}\right).
\end{align*}
If ${\partial^4u_0\over\partial t^4}\in L^2(0,T;H)$ and ${\partial^4\nabla_{y_i}\fu_i\over\partial t^4}\in L^2(0,T;H_i)$, then there is a constant $c$ such that

\begin{align*}
&\|\partial_t u^L_0-\partial_tu_0\|_{\tL(0,T;H)}+\sum_{i=1}^n\|\partial_t\nabla_{y_i}\fu_i^L-\partial_t\nabla_{y_i}\fu_i\|_{\tL(0,T;H_i)}\\
&\qquad+\|\curl u^L_{0}-\curl u_0\|_{\tL(0,T;H)}+\sum_{i=1}^n\|\scurl_{y_i}u^L_{i}-\scurl_{y_i}u_i\|_{\tL(0,T;H_i)}\\
&\quad\le c\left[(\Delta t)^2 \left\|{\partial^4u_0\over\partial t^4}\right\|_{L^2(0,T;H)}+(\Delta t)^2\sum_{i=1}^n\left\|{\partial^4\nabla_{y_i}\fu_i\over\partial t^4}\right\|_{H_i}+\left\|{\partial^2q^L_0\over\partial t^2}\right\|_{L^2(0,T;H)}\right.\\
&\qquad\left.+\sum_{i=1}^n\left\|{\partial^2\over\partial t^2}\nabla_{y_i}\fq_{i}^L\right\|_{L^2(0,T;H_i)}+\|q_0^L\|_{L^\infty(0,T;H)}+\sum_{i=1}^n\|\nabla_{y_i}\fq^L_i\|_{L^\infty(0,T;H_i)}\right]\\
&\qquad+c\left(\|\partial_t p^L_{0,1/2}\|_{H}+\sum_{i=1}^n\|\partial_t\nabla_{y_i} \fp^L_{i,1/2}\|_{H_i}+\|\curl p^L_{0,1/2}\|_H+\sum_{i=1}^n\|\scurl_{y_i}p^L_{i,1/2}\|_{H_i}\right)\\
&\qquad+\|\partial_tq_0^L\|_{\tL(0,T;H)}+\|\curl q_0^L\|_{\tL(0,T;H)}\\
&\qquad+\sum_{i=1}^n\left(\|\partial_t\nabla_{y_i}\fq_i^L\|_{\tL(0,T;H_i)}+\|\scurl_{y_i}q_i^L\|_{\tL(0,T;H_i)}\right).
\end{align*}
\end{proposition}
\bproof
We note that $\bu^L-\bu=\boldp^L+\boldq^L$. The conclusions follow from Lemma \ref{lem:fullydiscpest}.
\eproof

From this, we deduce
\bp
%\label{prop:semidiscconv}
If $\bu\in H^2(0,T;\bV)$, ${\partial^3u_0\over\partial t^3}\in L^2(0,T;H)$ and ${\partial^3\nabla_{y_i}\fu_i\over\partial t^3}\in$ $ L^2(0,T;H_i)$, and if we choose $\bu_0^L$ and $\bu_1^L$ such that
\[
\lim_{L\to 0}\|\partial_t p^L_{0,1/2}\|_{H}+\sum_{i=1}^n\|\partial_t\nabla_{y_i} \fp^L_{i,1/2}\|_{H_i}+\|\curl p^L_{0,1/2}\|_H+\sum_{i=1}^n\|\scurl_{y_i}p^L_{i,1/2}\|_{H_i}=0,
\]
then 
\begin{multline*}
\lim_{L\to\infty}\|\partial_t u^L_0-\partial_tu_0\|_{\tL(0,T;H)}+\sum_{i=1}^n\|\partial_t\nabla_{y_i}\fu_i^L-\partial_t\nabla_{y_i}\fu_i\|_{\tL(0,T;H_i)}\\
+\|\curl u^L_{0}-\curl u_0\|_{\tL(0,T;H)}+\sum_{i=1}^n\|\scurl_{y_i}u^L_{i}-\scurl_{y_i}u_i\|_{\tL(0,T;H_i)}=0.
\end{multline*}
\epr
\bproof
From the hypothesis and Lemma \ref{lem:estdtqL}, we have that
\[
\lim_{L\to\infty}\left\|{\partial^2\boldq^L\over\partial t^2}\right\|_{L^2(0,T;\bV)}=0.
\]
As $\bu \in C([0,T], \bV)$, from the proof of Proposition \ref{prop:semidiscerror} $\lim_{L\to\infty}\|\boldq^L\|_{L^\infty(0,T;\bV)}=0.$
We have that 
\[
\|\boldq^L\|_{\tL(0,T;\bV)}\le \|\boldq^L\|_{L^\infty(0,T;\bV)}
\]
so  
\[
\lim_{L\to\infty}\|\boldq^L\|_{\tilde{L}^\infty(0,T;\bV)}=0.
\]
Further, from \eqref{eq:eqforbwL}, we have that
\[
B(\partial_t\bw^L_{m+1/2}-\partial_t\bu_{m+1/2},\bv^L)+A(\partial_t\bw^L_{m+1/2}-\partial_t\bu_{m+1/2},\bv^L)=0
\]
for all $\bv^L\in\bV^L$. We thus have 
\[
\|\partial_t\bw^L_{m+1/2}-\partial_t\bu_{m+1/2}\|_\bV\le c\inf_{\bv^L\in \bV^L}\|\bv^L-\partial_t\bu_{m+1/2}\|_{\bV}.
\]
As $\partial_t\bu_{m+1/2}={\partial\bu\over\partial t}(\xi)$ for $\xi\in (0,T)$, we deduce that
\[
\|\partial_t\boldq^L\|_{\tL(0,T;\bV)}\le c\sup_{t\in (0,T)}\inf_{\bv^L\in \bV^L}\|\bv^L-{\partial\bu\over\partial t}(t)\|_{\bV}.
\]
As ${\partial\bu\over\partial t}\in C([0,T];\bV)$, a proof identical to that for $\|\boldq^L\|_{L^\infty(0,T;\bV)}$ in Proposition \ref{prop:semidiscerror} shows that
\[
\lim_{L\to\infty}\|\partial_t\boldq^L\|_{\tL(0,T;\bV)}=0.
\]
We thus get the conclusion. \eproof

\section{Regularity of the solution}
To derive an explicit error estimate for the full and sparse tensor product finite element approximating problems in the next section, 
we now establish the regularity of $u_0$ and $\nabla_{y_i}\fu_i$ with respect to $t$. 
%We make the following assumption on the smoothness of the matrix functions $a(x,\by)$ and $b(x,\by)$.
%
The function $\fu_i$ and $u_i$ can be written in terms of $u_0$ from the solution of the cell problems. Let $b^n(x,\by_n)=b(x,\by)$. Recursively, for all $i=0,\ldots,n$, let 
%\[
%\int_0^T \fu_i q''(t)dt=\left(\int_0^T(u_0+\sum_{j=1}^{i-1}\nabla_{y_j}\fu_j)_kq''(t)dt\right)w^k_i,
%\]
%where 
$w^k_i\in V_i$ be the solution of the cell problem
\be
\nabla_{y_i}\cdot(b^i(x,\by_i)(e^k+\nabla_{y_i} w^k_i))=0
\label{eq:cellbi}
\ee
where $e^k$ is the $k$th unit vector with every component equals 0, except the $k$th component which equals 1. 
%From an argument as above, we have 
%\[
%{\partial\over\partial t}\nabla_{y_i}\fu_i={\partial\over\partial t}(u_0+\sum_{j=1}^{i-1}\nabla_{y_j}\fu_j)_k\nabla_{y_i} w^k_i+G_i(x,\by_i)
%\]
%for a function $G_i(x,\by_i)\in L^2(D\times\bY_i)$. 
For $i=1,\ldots,n$, the positive definite matrix function $b^{i-1}(x)$ is defined as
\be
b^{i-1}_{pq}(x,\by_{i-1})=\int_{Y_i}b_{kl}^{i}(x,\by_{i})\left(\delta_{ql}+{\partial w^q_i\over\partial y_{il}}\right)\left(\delta_{pk}+{\partial w^p_i\over\partial y_{ik}}\right)dy_i;
\label{eq:bi}
\ee
 $b^0$ is the homogenized coefficient. 
Let $a^n=a$.  
%we then have, recursively,
%\[
%u_{i}=N_{i}^l \big((\curl u_0)_l+(\scurl_{y_1}u_1)_l+...+(\scurl_{y_{i-1}}u_{i-1})_l\big)
%\]
Let $N_i^k \in  W_i$ be the solution of  the cell problem
\be
\scurl_{y_i}(a^i(x,\by_i)(e^k+\scurl_{y_i}N^k_i))=0.
\label{eq:cellai}
\ee
%i.e. 
%\[
%    \int_D\int_{\bY_{i}} a^i(e^l+\scurl_{y_{i}}N^l_{i})\cdot\scurl_{y_{i}}v_{i} d\by_i=0
%\]
%for all $v_{i} \in W_i$. For $i=1,\ldots,n-1$, 
For $i=1,\ldots,n$, the positive definite coefficient $a^{i-1}$ is defined as
\be
a^{i-1}_{pq}(x,\by_{i-1})=\int_{Y_{i}}a^i_{kl}(x,\by_i)\left(\delta_{ql}+(\scurl_{y_{i}}N_{i}^q)_l\right)\left(\delta_{pk}+(\scurl_{y_i}N_i^p)_k\right)dy_{i},
%=\int_{Y_n}a\left(I+\scurl_{y_{i+1}}N_{i+1}\right)\left(I+\scurl_{y_i}N_{i+1}\right)dy_{i+1};
\label{eq:ai}
\ee
$a^0$ is the homogenized coefficient.
%Continuing this process, we finally get the homogenized coefficient $a^0(x)$ as
%\be
%a^0_{pq}(x)=\int_{Y_{1}}a^1_{pk}\left(\delta_{kq}+(\scurl_{y_{1}}N_{1}^q)_k\right).
%\left(I+\scurl_{y_1}N_1\right)dy_{1}.
%\label{eq:a0}
%\ee
The homogenized equation is
\[
\int_0^T\int_D b^0u_0\cdot v_0q''(t)dxdt+\int_0^T\int_D a^0\curl u_0\cdot\curl v_0 q(t) dxdt=\int_0^T\int_D f(t,x)\cdot v_0(x) q(t)dxdt
\]
for all $q\in {\cal D}(0,T)$ and $v_0\in W$, 
i.e. 
\beq
b^0(x) \frac{\partial^2u_0}{\partial t^2}(t,x)+\curl(a^0(x)\curl u_0(t,x))=f(t,x).
\label{eq:homogenizedeq}
\eeq
The solution $\bu$ is written in terms of $u_0$ as
\be
u_i=N_i^{r_{i-1}}(\delta_{r_{i-1}r_{i-2}}+(\scurl_{y_{i-1}}N^{r_{i-2}}_{i-1})_{r_{i-1}})(\delta_{r_{i-2}r_{i-3}}+(\curl_{y_{i-2}}N^{r_{i-3}}_{i-2})_{r_{i-2}})\ldots(\delta_{r_1r_0}+(\curl_{y_1}N^{r_0}_1)_{r_1})(\curl u_0)_{r_0},
\label{eq:ui}
\ee
and
\begin{eqnarray}
\nabla_{y_i}\fu_i=u_{0r_0}(\delta_{r_0r_1}+{\partial w_1^{r_0}\over\partial y_{1r_1}})(\delta_{r_1r_2}+{\partial w_2^{r_1}\over\partial  y_{2r_2}})\ldots\nabla_{y_i}w^{r_{i-1}}_i-g_{1r_0}(\delta_{r_0r_1}+{\partial w_1^{r_0}\over\partial y_{1r_1}})(\delta_{r_1r_2}+{\partial w_2^{r_1}\over\partial  y_{2r_2}})\ldots\nabla_{y_i}w^{r_{i-1}}_it\nonumber\\
-g_{0r_0}(\delta_{r_0r_1}+{\partial w_1^{r_0}\over\partial y_{1r_1}})(\delta_{r_1r_2}+{\partial w_2^{r_1}\over\partial  y_{2r_2}})\ldots\nabla_{y_i}w^{r_{i-1}}_i.
\label{eq:fui}
\end{eqnarray}
We refer to \cite{Mwt} for detailed derivation. We make the following assumption on the smoothness of the matrix functions $a(x,\by)$ and $b(x,\by)$.

\begin{assumption}
\label{assum:regularityab}
The matrix functions $a$ and $b$ belong to $C^1(\bar D,C^2(\bar Y_1,\ldots,C^2(\bar Y_n)\ldots))^{d\times d}$.
\end{assumption}
With this assumption, we  have
\begin{proposition}
\label{prop:regularityNomega}
Under Assumption \ref{assum:regularityab}, for all $i,r=1,\ldots,d$, $\scurl_{y_i}N_i^r\in C^1(\bar D,C^2(\bar Y_1,\ldots,C^2(\bar Y_{i-1},H^2(Y_i))\ldots))$ and $w_i^r\in C^1(\bar D,C^2(\bar Y_1,\ldots,C^2(\bar Y_{i-1},H^3(Y_i))\ldots))$.
\end{proposition}
We refer to  \cite{Tiep1} for a proof of this proposition.
We have the following regularity results for the solution $u_0$ of the homogenized equation \eqref{eq:homogenizedeq}. 

\begin{proposition}
Under Assumption \ref{assum:regularityab}, assume 
\beq
\begin{cases}
f \in H^2(0,T; H),&\\
g_1 \in W,&\\
(b^0)^{-1}[f(0)-\curl (a^0(x)\curl g_0)] \in W,&\\
(b^0)^{-1}[\frac{\partial f}{\partial t}(0)-\curl (a^0(x)\curl g_1)] \in H,&
\end{cases}
\label{eq:Compatibility1}
\eeq
then 
\beq
\frac{\partial^2 u_0}{\partial t^2}\in L^\infty(0,T; W),\ \ \frac{\partial^3 u_0}{\partial t^3}\in L^\infty(0,T; H),\ \ \mbox{and}\ \ {\partial^3\over\partial t^3}\nabla_{y_i}\fu_i\in L^\infty(0,T;L^2(D\times\bY)).
\label{eq:reg1}
\eeq
Further, if
\beq
\begin{cases}
f \in H^3(0,T; H),&\\
g_1 \in W,&\\
(b_0)^{-1}[f(0)-\curl (a^0(x)\curl g_0)] \in W,&\\
(b_0)^{-1}[\frac{\partial f}{\partial t}(0)-\curl (a^0(x)\curl g_1)] \in W,&\\
(b_0)^{-1}[\frac{\partial^2 f}{\partial t^2}(0)-\curl (a^0(x)\curl ((b_0)^{-1}(f(0)-\curl (a^0(x)\curl g_0)))) \in H,&
\end{cases}
\label{eq:Compatibility2}
\eeq
then 
\beq
\frac{\partial^3 u_0}{\partial t^3}\in L^\infty(0,T; W),\ \ \frac{\partial^4 u_0}{\partial t^4}\in L^\infty(0,T; H),\ \ \mbox{and}\ {\partial^4\over\partial t^4}\nabla_{y_i}\fu_i\in L^\infty(0,T;L^2(D\times\bY)).  
\label{eq:reg2}
\eeq
\end{proposition}

\bproof
We use the regularity theory of general hyperbolic equations (see, e.g., Wloka \cite{Wloka}, Chapter 5). 
From \eqref{eq:Compatibility1} we have that
\be
b^0\frac{\partial ^2}{\partial t ^2}\left (\frac{\partial u_0}{\partial t }\right)+\curl \left(a^0\curl \frac{\partial u_0}{\partial t }\right)=\frac{\partial f}{\partial t }
\label{eq:e1}
\ee
with compatibility initial conditions
\[
\frac{\partial u_0}{\partial t }(0)=g_1 \in W,~~\frac{\partial }{\partial t }\frac{\partial u_0}{\partial t }(0)=(b^0)^{-1}[f(0)-\curl (a^0\curl g_0)] \in W
\]
and
\be
\label{eq:e2}
b^0\frac{\partial ^2}{\partial t ^2}\left(\frac{\partial^2 u_0}{\partial t^2 }\right)+\curl \left(a^0\curl \frac{\partial^2 u_0}{\partial t^2 }\right)=\frac{\partial^2 f}{\partial t^2 },
\ee
with compatibility initial conditions
\[
\frac{\partial^2 u_0}{\partial t^2 }(0)=(b^0)^{-1}[f(0)-\curl \left(a^0\curl g_0\right)] \in W
\ \ 
\mbox{and}
\ \ 
\frac{\partial }{\partial t }\frac{\partial^2 u_0}{\partial t^2 }(0)=(b^0)^{-1}[\frac{\partial f}{\partial t}(0)-\curl (a^0\curl g_1)] \in  H.
\]
We thus deduce that
\[
{\partial^2u_0\over\partial t^2}\in L^\infty(0,T;W)\ \ \mbox{and}\ \ {\partial^3u_0\over\partial t^3}\in L^\infty(0,T;H).
\]
From \eqref{eq:fui} and  Proposition \ref{prop:regularityNomega}, we deduce that
\[
{\partial^3\over\partial t^3}\nabla_{y_i}\fu_i\in L^\infty(0,T;L^2(D\times\bY)).
\]
Similarly, we deduce regularity \eqref{eq:reg2} from \eqref{eq:Compatibility2}. \eproof

To derive explicitly the rate of convergence for the full and sparse tensor finite element approximations in the next section, we define the following regularity spaces. 
For $i=1,\ldots,n$, 
let $\bar{\cal H}_i$ be the space of functions belonging to $L^2(D\times Y_1\times\ldots\times Y_{i-1},H^1_\#(\curl,Y_i))$, 
%such that for all $k=1,2,3$
%\[
%{\partial w\over\partial x_k}\in L^2(D\times Y_1\times\ldots\times Y_{i-1},\tilde H_\#(\curl,Y_i)),
%\]
%and for all $j=1,\ldots,i-1$ and $k=1,2,3$
%\[
%{\partial w\over\partial (y_j)_k}\in L^2(D\times Y_1\times\ldots\times Y_{i-1},\tilde H_\#(\curl,Y_i)).
%\]
%In other words, for all $w\in {\cal H}_i$, $w$ belongs to 
%$L^2(D\times Y_1\times\ldots\times Y_{i-1},H^1_\#(\curl,Y_i))\bigcup 
$L^2(Y_1\times\ldots\times Y_{i-1},H^1(D,\tilde H_\#(\curl,Y_i)))$ and $L^2(D\times\prod_{k<i,k\ne j}Y_k,H^1_\#(Y_j,\tilde H_\#(\curl,Y_i)))$ for $j=1,\ldots,i-1$. 
For $0<s<1$, by interpolation, we define the space $\bar{\cal H}^s_i$ which 
%\footnote{We need the regularity of the solutions of the cell problems to deduce the homogenization error in Section \ref{sec:homerror} ahead so we consider only the case of $H^1(Y_j)$ regularity with respect to $y_j$ for $j=1,\ldots,i-1$.} 
consists of functions $w$ that belongs to $L^2(D\times Y_1\times\ldots\times Y_{i-1},H^s_\#(\curl,Y_i))$, $L^2(Y_1\times\ldots\times Y_{i-1},H^s(D,\tilde H_\#(\curl,Y_i)))$ and $L^2(D\times\prod_{k<i,k\ne j}Y_k,H^s_\#(Y_j,\tilde H_\#(\curl,Y_i)))$.  We equip $\bar{\cal H}_i^s$ with the norm
\beqas
\|w\|_{\bar{\cal H}_i^s}=\|w\|_{L^2(D\times Y_1\times\ldots\times Y_{i-1},H^s_\#(\curl,Y_i))}+\|w\|_{L^2(Y_1\times\ldots\times Y_{i-1},H^s(D,\tilde H_\#(\curl,Y_i)))}+\\
\sum_{j=1}^{i-1}\|w\|_{L^2(D\times\prod_{k<i,k\ne j}Y_k,H^s_\#(Y_j,\tilde H_\#(\curl,Y_i)))}.
\eeqas
We define $\bar{\frak H}_i^s$ as the space of functions  $w\in L^2(D\times Y_1\times\ldots\times Y_{i-1},H^{1+s}_\#(Y_i))$ such that $w\in L^2(Y_1\times\ldots\times Y_{i-1},H^s(D,H^1_\#(Y_i)))$ and for all $j=1,\ldots,i-1$, $ w\in L^2(D\times\prod_{k<i,k\ne j}Y_k,H^s_\#(Y_{j},H^1_\#(Y_i)))$. We equip this space with the norm
\beqas
\|w\|_{\bar{\frak H}_i^s}=\|w\|_{L^2(D\times Y_1\times\ldots\times Y_{i-1},H^{1+s}_\#(Y_i))}+\|w\|_{L^2(Y_1\times\ldots\times Y_{i-1},H^s(D,H^1_\#(Y_i)))}+\\
\displaystyle \sum_{j=1}^{i-1}\|w\|_{L^2(D\times\prod_{k<i,k\ne j}Y_k,H^s_\#(Y_j,H^1_\#(Y_i)))}.
\eeqas
We define the regularity space $\bar{\bH^s}$ as
\[
\bar{\bH^s}=H^s(\curl,D)\times\bar{\cal H}_1^s\times\ldots\bar{\cal H}_n^s\times\bar{\frak H}_1^s\times\ldots\times\bar{\frak H}_n^s.
\]
We define $\hat{\cal H}_i$ as the space of functions $w\in L^2(D\times Y_1\times\ldots\times Y_{i-1},H^1_\#(\curl,Y_i))$ which are periodic with respect to $y_j$ with the period being $Y_j$ ($j=1,\ldots,i-1$) such that for any $\alpha_0,\alpha_1,\ldots\alpha_{i-1}\in \IR^d$ with $|\alpha_k|\le 1$ for $k=0,\ldots,i-1$, 
\[
{\partial^{|\alpha_0|+|\alpha_1|+\ldots+|\alpha_{i-1}|}\over\partial x^{\alpha_0}\partial y_1^{\alpha_1}\ldots\partial y_{i-1}^{\alpha_{i-1}}}w\in L^2(D\times Y_1\times\ldots\times Y_{i-1},H^1_\#(\curl,Y_i)).
\]
We equip $\hat{\cal H}_i$ with the norm
\[
\|w\|_{\hat{\cal H}_i}=\sum_{\alpha_j\in \IR^d,|\alpha_j|\le 1\atop 0\le j\le i-1}\left\|{\partial^{|\alpha_0|+|\alpha_1|+\ldots+|\alpha_{i-1}|}\over\partial x^{\alpha_0}\partial y_1^{\alpha_1}\ldots\partial y_{i-1}^{\alpha_{i-1}}}w\right\|_{ L^2(D\times Y_1\times\ldots\times Y_{i-1},H^1_\#(\curl,Y_i))}.
\]
We can write $\hat{\cal H}_i$ as $H^1(D,H^1_\#(Y_1,\ldots,H^1_\#(Y_{i-1},H^1_\#(\curl,Y_i))))$. 

By interpolation, we define $\hat{\cal H}_i^s=H^s(D,H^s_\#(Y_1,\ldots,H^s_\#(Y_{i-1},H^s_\#(\curl,Y_i))\ldots))$ for $0<s<1$.

We define $\hat{\frak H}_i$ as the space of functions $w\in L^2(D\times Y_1\times\ldots\times Y_{i-1},H^2_\#(Y_i))$ that are periodic with respect to $y_j$ with the period being $Y_j$ for $j=1,\ldots,i-1$ such that $\alpha_0,\alpha_1,\ldots\alpha_{i-1}\in \IR^d$ with $|\alpha_k|\le 1$ for $k=0,\ldots,i-1$, 
\[
{\partial^{|\alpha_0|+|\alpha_1|+\ldots+|\alpha_{i-1}|}\over\partial x^{\alpha_0}\partial y_1^{\alpha_1}\ldots\partial y_{i-1}^{\alpha_{i-1}}}w\in L^2(D\times Y_1\times\ldots\times Y_{i-1},H^2_\#(Y_i)).
\]
The space $\hat{\frak H}_i$ is equipped with the norm
\[
\|w\|_{\hat{\frak H}_i}=\sum_{\alpha_j\in \IR^d,|\alpha_j|\le 1\atop 0\le j\le i-1}\left\|{\partial^{|\alpha_0|+|\alpha_1|+\ldots+|\alpha_{i-1}|}\over\partial x^{\alpha_0}\partial y_1^{\alpha_1}\ldots\partial y_{i-1}^{\alpha_{i-1}}}w\right\|_{ L^2(D\times Y_1\times\ldots\times Y_{i-1},H^2_\#(Y_i))}.
\]
We can write $\hat{\frak H}_i$ as $H^1(D,H^1_\#(Y_1,\ldots,H^1_\#(Y_{i-1},H^2_\#(Y_i))))$. By interpolation,  we define the space\\
 $\hat{\frak H}_i^s:=H^s(D,H^s(Y_1,\ldots,H^s(Y_{i-1},H^{1+s}_\#(Y_i))))$. The regularity space $\hat{\bH^s}$ is defined as
\[
\hat{\bH^s}=H^s(\curl,D)\times\hat{\cal H}_1^s\times\ldots\hat{\cal H}_n^s\times\hat{\frak H}_1^s\times\ldots\times\hat{\frak H}_n^s.
\]

For the regularity of $u_0$, we have the following result.
\bp\label{prop:regularityu0}
Under Assumption \ref{assum:regularityab}, if $D$ is a Lipschitz polygonal domain, $f\in H^1(0,T;H)$, $g_0\in H^1(\curl,D)$ and $g_1\in W$, ${\rm div}f\in L^\infty(0,T;L^2(D))$, ${\rm div}(b^0g_0)\in L^2(D)$ and ${\rm div}(b^0g_1)\in L^2(D)$,  there is a constant $s\in(0,1]$ such that $u_0\in L^\infty(0,T;H^s(\curl,D))$.
\epr
\bproof
Using Proposition \ref{prop:regularityNomega}, equations \eqref{eq:ai} and \eqref{eq:bi}, we have that $a^0, b^0\in C^1(\bar D)^{d\times d}$.
As $f\in H^1(0,T;H)$ and $g_0\in H^1(\curl,D)$, we have that $(b^0)^{-1}[f-\curl(a^0\curl g^0)]\in H$. The compatibility initial conditions hold so that ${\partial^2u_0\over\partial t^2}\in L^\infty(0,T;H)$. Thus
\[
\curl(a^0\curl u_0)=f-b^0{\partial^2u_0\over\partial t^2}\in L^\infty(0,T;H).
\]
Let $U(t)=a^0\curl u_0(t)$. As ${\rm div}((a^0)^{-1}U(t))=0$ and $(a^0)^{-1}U(t)\cdot\nu=0$, there is a constant $c$ and a constant $s\in (0,1]$ which depend on $a^0$ and the domain $D$ so that
\[
\|U(t)\|_{H^s(D)^3}\le c(\|\curl U(t)\|_{L^2(D)^3}+\|U(t)\|_{L^2(D)^3})
\]
so $U\in L^\infty(0,T;H^s(D)^3)$. As $\curl u_0(t)=(a^0)^{-1}U(t)$ and $(a^0)^{-1}\in C^1(\bar D)^{d\times d}$, $\curl u_0\in L^\infty(0,T;H^s(D))$. 
We note that
\[
{\rm div}\left(b^0{\partial^2u_0\over\partial t^2}\right)={\rm div} f,
\]
so 
\[
{\rm div}(b^0u_0(t))=\int_0^t\int_0^s{\rm div}f(r)drds+t{\rm div}(b^0g_1)+{\rm div}(b^0 g_0)\in L^\infty(0,T;L^2(D)).
\]
From Theorem 4.1 of Hiptmair \cite{Hiptmair2002}, we deduce that there is a constant $s\in (0,1]$ (we take it as the same constant as above), so that 
\[
\|u_0(t)\|_{H^s(D)^3}\le c(\|u_0(t)\|_{H(\curl,D)}+\|{\rm div}(b^0u_0(t))\|_{L^2(D)^3}).
\]
Thus $u_0\in L^\infty(0,T;H^s(\curl,D))$.   
\eproof

Similarly, we can deduce the regularity for ${\partial^2u_0\over\partial t^2}$.
\bp\label{prop:regularityd2u0dt2}
Under Assumption \ref{assum:regularityab}, if $D$ is a Lipschitz polygonal domain, if the compatibility conditions \eqref{eq:Compatibility2} hold, and if ${\rm div}f\in L^\infty(0,T;L^2(D))$, then  there is a constant $s\in (0,1]$ such that ${\partial^2u_0\over\partial t^2}\in L^\infty(0,T;H^s(\curl,D))$.
\epr
\bproof
From equation \eqref{eq:e2}, we have
\[
\curl\left(a^0\curl{\partial^2u_0\over\partial t^2}\right)={\partial^2 f\over\partial t^2}-b^0{\partial^4u_0\over\partial t^4}\in L^\infty(0,T;H)
\]
as ${\partial^4u_0\over\partial t^4}\in L^\infty(0,T;H)$ due to \eqref{eq:Compatibility2}. Following a similar argument as in the proof of Proposition \ref{prop:regularityu0} we deduce that $\curl{\partial^2u_0\over\partial t^2}\in L^\infty(0,T;H^s(D)^3)$. We note that 
\[
{\rm div}b^0{\partial^2 u_0\over\partial t^2}={\rm div} f\in L^\infty(0,T;L^2(D)).
\]
From Theorem 4.1 of \cite{Hiptmair2002}, we deduce that ${\partial^2u_0\over\partial t^2}\in L^\infty(0,T;H^s(D)^3)$.  
\eproof

From these we deduce 
\bp 
Under Assumption \ref{assum:regularityab}, and the hypothesis of Proposition \ref{prop:regularityu0}, there is a constant $s\in (0,1]$ so that $\bu\in L^\infty(0,T;\hat{\bH^s})$.
\epr
\bproof
From Proposition \ref{prop:regularityNomega}, we have that $N_i^r$ and $\curl N_i^r$ belong to $C^1(\bar D,C^2(\bar Y_1,\ldots,C^2(\bar Y_{i-1},H^2_\#(Y_i))\ldots))$. Together with $u_0\in L^\infty(0,T;H^s(\curl,D))$, this implies $u_i\in L^\infty(0,T;\hat{\cal H}_i^s)$. Similarly, we have $\fu_i\in L^\infty(0,T;\hat{\frak H}_i^s)$. 
\eproof

Similarly, we have:
\bp \label{prop:boldu_tt}
Under Assumption \ref{assum:regularityab}, and the hypothesis of Proposition \ref{prop:regularityd2u0dt2}, there is a constant $s\in (0,1]$ so that ${\partial^2\bu\over\partial t^2}\in L^\infty(0,T;\hat{\bH^s})$.
\epr

\begin{remark}\label{rem:regularitybud2tbu} We have
\[
\scurl {\partial u_0(t)\over\partial t}=\int_0^t\curl{\partial^2u_0\over\partial t^2}(s)ds+\curl g_1,\ \ \mbox{and}\ \ 
{\partial u_0(t)\over\partial t}=\int_0^t{\partial^2u_0\over\partial t^2}(s)ds+g_1,
\]
\[
\curl u_0(t)=\int_0^t\int_0^s\curl{\partial^2u_0\over\partial t^2}(r)drds+t\curl g_1+\curl g_0,\ \ \mbox{and}\ \ 
u_0(t)=\int_0^t\int_0^s{\partial^2u_0\over\partial t^2}(r)drds+tg_1+g_0.
\]
Thus with the hypothesis of Proposition \ref{prop:regularityd2u0dt2}, together with $g_0\in H^s(\curl,D)$ and $g_1\in H^s(\curl,D)$, we deduce that ${\partial u_0\over\partial t}\in L^\infty(0,T;H^s(\curl,D))$ and $u_0\in L^\infty(0,T;H^s(\curl,D))$. This implies also that $\bu\in L^\infty(0,T;\hat{\bH^s})$. 
\end{remark}

\section{Full and sparse tensor product approximations}
We consider the approximations of problem \eqref{eq:msprob} using the full and sparse tensor product FE. 
We assume that the domain  $D$ is a polygon in $\IR^3$. Let ${\cal T}^l$ ($l=0,1,\ldots$) be the sets of simplices in $D$ with mesh size $h_l=O(2^{-l})$ which are determined recursively where ${\cal T}^{l+1}$ is obtained from ${\cal T}^l$ by dividing each simplex in ${\cal T}^l$ into 
%to 4 congruent triangles for $d=2$ or 
8 tedrahedra. For a tedrahedron $T\in {\cal T}^l$, we consider the edge finite element space
\[
R(T)=\{v:\ \ v=\alpha+\beta\times x,\ \ \alpha,\beta\in\IR^3\}.
\]
When $D$ is a polygon in $\IR^2$, ${\cal T}^{l+1}$ is obtained from ${\cal T}^l$ by dividing each simplex in ${\cal T}^l$ into 4 congruent triangles. For each triangle $T\in {\cal T}^l$, we consider the edge finite element space
\[
R(T)=\left\{v:\ \ v=\left(
\begin{array}{c}
\alpha_1\\
\alpha_2
\end{array}\right)
+
\beta\left(\begin{array}{c}
x_2\\
-x_1
\end{array}\right)
\right\}
\]
where $\alpha_1,\alpha_2$ and $\beta$ are constants. 
Alternatively, when $D$ is partitioned into cubic meshes, we can use edge finite element on cubic mesh instead (see \cite{Monk}). For each simplex $T\in {\cal T}^l$, we denote by ${\cal P}_1(T)$ the set of linear polynomials in ${\cal T}$. In the following, we only present the analysis for the three dimensional case as the two dimensional case is similar.

We define the finite element spaces
\beqas
W^l=\{v\in H_0(\curl,D),\ v|_T\in R(T)\ \forall\,T\in {\cal T}^l\},\\
V^l=\{v\in H^1(D),\ v|_T\in {\cal P}_1(T)\ \forall\,T\in {\cal T}^l\}.
\eeqas
For the cube $Y$, we consider a hierarchy of simplices ${\cal T}^l_\#$ that are distributed periodically. We consider the space of functions 
\beqas
W^l_\#=\{v\in H_\#(\curl,Y),\ v|_T\in R(T)\ \forall\,T\in {\cal T}^l_\#\}
\eeqas
and 
\[
V_\#^l=\{v\in H^1_\#(Y),\ v|_T\in {\cal P}_1(T)\ \forall\,T\in {\cal T}^l_\#\}.
\]
 We then have the following standard estimates (see Monk \cite{Monk} and Ciarlet \cite{Ciarlet})
\[
\inf_{v_l\in W^l}\|v-v_l\|_{H(\curl,D)}\le ch_l^s(\|v\|_{H^s(D)^d}+\|\curl v\|_{H^s(D)^d})
\]
$\forall\,v\in H_0(\curl, D)\bigcap H^s(\curl,D)$;
\[
\inf_{v_l\in W_\#^l}\|v-v_l\|_{H_\#(\curl,Y)}\le ch_l^s(\|v\|_{H^s(Y)^d}+\|\curl v\|_{H^s(Y)^d})
\]
$\forall\,v\in H_\#(\curl, Y)\bigcap H^s(\curl, Y)$;
\[
\inf_{v_l\in V^l}\|v-v_l\|_{L^2(D)}\le ch_l^s\|v\|_{H^s(D)}
\]
$\forall\,v\in H^s(D)$;
\[
\inf_{v_l\in V^l_\#}\|v-v_l\|_{L^2(Y)}\le ch_l^s\|v\|_{H^{s}(Y)}
\]
$\forall\,v\in H^s_\#(Y)$;
 and
\[
\inf_{v_l\in V^l_\#}\|v-v_l\|_{H^1_\#(Y)}\le ch_l^s\|v\|_{H^{1+s}(Y)}
\]
$\forall\,v\in H^1_\#(Y)\bigcap H^{1+s}(Y)$. 
\subsection{Full tensor product finite elements} \label{sec:fulltensor}
As $L^2(D\times\bY_{i-1},\tilde H_\#(\curl,Y_i))\cong L^2(D)\otimes L^2(Y_1)\otimes\ldots\otimes L^2(Y_{i-1})\otimes \tilde H_\#(\curl,Y_i)$ we use the tensor product finite element space
\[
\bar W_i^l=V^l\otimes \underbrace{V^l_\#\otimes\ldots\otimes V^l_\#}_{i-1\mbox{ times }}\otimes W^l_\#
\]
to approximate $u_i$. Similarly, as $\frak u_i\in L^2(D\times\bY_{i-1},H^1_\#(Y))$, we use the finite element space
\[
\bar V_i^l=V^l\otimes \underbrace{V^l_\#\otimes\ldots\otimes V^l_\#}_{i\mbox{ times }}
\]
to approximate $\frak u_i$. We define the space
\[
\bar{\bV}^L=W^L\otimes\bar W_1^L\otimes\ldots\otimes\bar W_n^L\otimes \bar V_1^L\otimes\ldots\otimes\bar V_n^L.
\]
The spatially semidiscrete  full tensor product finite element approximating problem is: Find $\bbu^L(t)\in\bbV^L$ so that for all $\bbv^L\in \bbV^L$:
\begin{align}
&\int_D\int_\bY b(x,\by)\left[\left({\partial^2\bau^L_0\over \partial t^2}(t)+\sum_{i=1}^n{\partial^2\over \partial t^2}\nabla_{y_i} \bafu^L_i(t)\right)\cdot\left(\bav_0^L+\sum_{i=1}^n\nabla_{y_i} \bafv^L_i\right)\right.d\by dx\nonumber\\
&\qquad+\left.a(x,\by)\left(\curl \bau^L_0(t)+\sum_{i=1}^n\scurl_{y_i}\bau^L_i(t)\right)\cdot\left(\curl \bav^L_0+\sum_{i=1}^n\scurl_{y_i}\bav^L_i\right)\right]d\by dx\nonumber\\
&\quad=\int_D f(t,x)\cdot \bav^L_0(x)dx
\label{eq:fullsemidiscprob}
\end{align}
 for all $\bbv^L=(\bav_0^{L},\{\bav_i^{L}\},\{\bafv_i^{L}\})\in \bbV^L$.
 
% We first state the following regularity for $\bu$
%\bp 
%Under Assumption \ref{assum:regularityab}, and the hypothesis of Proposition \ref{prop:regularityu0}, there is a constant $s\in (0,1]$ so that $\bu\in L^\infty(0,T;\hat{\bH^s})$.
%\epr
%\bproof
%From Proposition \ref{prop:regularityNomega1}, we have that $\chi_i^r$ and $\curl \chi_i^r$ belong to $$C^1(\bar D,C^2(\bar Y_1,\ldots,C^2(\bar Y_2,H^2(Y_i))\ldots))^3.$$ Together with $u_0\in L^\infty(0,T;H^s(\curl,D))$, this implies $u_i\in L^\infty(0,T;\hat{\cal H}_i^s)$. Similarly, we have $\fu_i\in L^\infty(0,T;\hat{\frak H}_i^s)$. 
%\eproof
%
%By the same argument, we deduce
%\bp 
%\label{prop:boldu_tt}
%Under Assumption \ref{assum:regularityab}, and the hypothesis of Proposition \ref{prop:regularityd2u0dt2}, there is a constant $s\in (0,1]$ so that ${\partial^2\bu\over\partial t^2}\in L^\infty(0,T;\hat{\bH^s})$.
%\epr
To deduce an error estimate for the full tensor product approximations of \eqref{eq:msprob}, we note the following approximations
\begin{lemma}
For $w\in \bar{\cal H}_i^s$,
\[
\inf_{w^l\in W_i^l}\|w-w^l\|_{L^2(D\times Y_1\times\ldots\times Y_{i-1},\tilde H_\#(\curl,Y_i))}\le ch_l^s\|w\|_{{\cal H}_i^s}.
\]
For $w\in \bar{\frak H}_i^s$, 
\[
\inf_{w^l\in V_i^l}\|w-w^l\|_{L^2(D\times Y_1\times\ldots\times Y_{i-1},H^1_\#(Y_i))}\le ch_l^s\|w\|_{\frak H_i^s}.
\]
\end{lemma}
The proofs of these results are similar to those for full tensor product finite elements in \cite{HSelliptic} and \cite{BungartzGriebel}, using orthogonal projection. We refer to \cite{HSelliptic} and \cite{BungartzGriebel} for details. From this we deduce that for $\bw\in\bar{\bH^s}$
\[
\inf_{\bw^L\in\bar\bV^L}\|\bw-\bw^L\|_{\bV}\le ch_L^s\|\bw\|_{\bar{\bH^s}}.
\]
We then have the following result for the spatially semidiscrete  approximation.
\begin{proposition}
\label{prop:fullsemidiscerror}
Assume that condition (\ref{eq:Compatibility2}) and Assumption \ref{assum:regularityab} hold, $D$ is a Lipschitz polygonal domain, ${\rm div}f\in L^\infty(0,T;L^2(D))$ and $g_0, g_1$ belong to $H^s(\curl,D)$. If  $g_0^L$ and $g_1^L$ are chosen so that
\be
\|g_0^L-g_0\|_W \leq ch_L^s \text{ and }     \|g_1^L-g_1\|_{H} \leq ch_L^s, 
\label{eq:fullsemidiscinitials}
\ee
where $s\in (0,1]$ is the constant in Proposition \ref{prop:regularityd2u0dt2}.
Then
\beqas
&&  \left\|\frac{\partial (\bau_0^L-u_0)}{\partial t}\right\|_{L^\infty(0,T;H)}+\sum_{i=1}^n\left\|\nabla_{y_i}\frac{\partial (\bafu^L_i-\fu_i)}{\partial t}\right\|_{L^\infty(0,T; H_i)}\\
&&+\|\scurl (\bau_0^L-u_0)\|_{L^\infty(0,T;H)}+\sum_{i=1}^n\|\scurl_{y_i}(\bau_i^L-u_i)\|_{L^\infty(0,T;H_i)}\le ch_L^s. 
\eeqas
\end{proposition}
\bproof From Proposition \ref{prop:boldu_tt} and Remark \ref{rem:regularitybud2tbu}, we deduce that $\bu\in L^\infty(0,T;\bar{\bH^s})$, ${\partial\bu\over\partial t}\in L^\infty(0,T;\bar{\bH^s})$ and ${\partial^2\bu\over\partial t^2}\in L^\infty(0,T;\bar{\bH^s})$. From Lemmas
\ref{lem:estqL} and \ref{lem:estdtqL}, we have
\be
\|\boldq^L\|_{L^\infty(0,T;\bV)}\le ch_L^s,\ \ \left\|{\partial\boldq^L\over\partial t}\right\|_{L^\infty(0,T;\bV)}\le ch_L^s,\ \ \mbox{and}\ \ \left\|{\partial^2\boldq^L\over\partial t^2}\right\|_{L^2(0,T;\bV)}\le ch_L^s.
\label{eq:estqLfull}
\ee
These together with
\begin{align*}
{\partial p_0^L\over\partial t}(0)&={\partial\over\partial t}(u_0^L(0)-u(0))-{\partial q_0^L\over\partial t},\\
\nabla_{y_i}{\partial p_i^L\over\partial t}(0)&=\nabla_{y_i}{\partial\over\partial t}(\fu_i^L(0)-\fu_i(0))-\nabla_{y_i}{\partial\fq_i^L\over\partial t}(0),\\
\curl p_0^L(0)&=\curl (u_0^L(0)-u_0(0))-\curl q_0^L(0),
\end{align*}
and \eqref{eq:fullsemidiscinitials}, we have that
\[
\left\|{\partial p_0^L\over\partial t}(0)+\sum_{i=1}^n\nabla_{y_i}{\partial p_i^L\over\partial t}(0)\right\|_{H_n}\le ch_L^s,
\]
and 
\[
\|\curl p_0^L(0)\|_{H}\le ch_L^s.
\]
Thus the right hand side of \eqref{eq:semidiscerrorest} is not more than $ch_L^s$. We thus get the conclusion.
\eproof

The fully  discrete problem now becomes:
For $m=1,...,M$ find $$\bbu^L_m=(\bau_{0,m}^L,\bau_{1,m}^L,...,\bau_{n,m}^L, \bar{\fu}_{1,m}^L ,...,\bar{\fu}_{n,m}^L) \in \bbV^L$$ such that for $m=1,...,M-1$
\begin{align}
&\int_D\int_\bY \left[b(x,\by)\left(\partial^2_t\bau^L_{0,m}+\sum_{i=1}^n\nabla_{y_i} \partial^2_t\bafu^L_{i,m}\right)\cdot\left(\bav_0^L+\sum_{i=1}^n\nabla_{y_i} \bafv^L_i\right)\right.\nonumber\\
&\qquad\left.+a(x,\by)\left(\curl \bau^L_{0 ,m , 1/4}+\sum_{i=1}^n\scurl_{y_i}\bau^L_{i ,m , 1/4}\right)\cdot\left(\curl \bav^L_0+\sum_{i=1}^n\scurl_{y_i}\bav^L_i\right)\right]d\by dx\nonumber\\
&\quad=\int_D f_{m,1/4}(x)\cdot \bav^L_0(x)dx
\label{eq:fullfullydiscprob}
\end{align}
for all $\bbv^L=(\bav_0^L, \bav_1^L,..., \bav^L_n, \bafv_1^L,...,\bafv_n^L) \in \bar\bV^L.$

\begin{proposition}
\label{prop:fullfullydiscerror}
Assume that condition (\ref{eq:Compatibility2}) and Assumption \ref{assum:regularityab} hold, $D$ is a Lipschitz polygonal domain, ${\rm div}f\in L^\infty(0,T;L^2(D))$ and $g_0, g_1$ belong to $H^s(\curl,D)$ where $s\in (0,1]$ is the constant in Proposition \ref{prop:regularityd2u0dt2}.
If the initial value $\bar\bu_1^L$ is chosen so that
\[
\|\partial_tp^L_{0,1/2}\|_{H}+\sum_{i=1}^n\|\partial_t\nabla_{y_i}\fp^L_{i,1/2}\|_{H_i}+\|\curl p^L_{0,1/2}\|_{H}+\sum_{i=1}^n\|\scurl_{y_i}p^L_{i,1/2}\|_{H_i}\le c((\Delta t)^2+h_L^s),
\] 
then 
\begin{multline*}
\|\partial_t\bau_0^L-\partial_tu_0\|_{\tL(0,T;H)}+\sum_{i=1}^n\left\|\partial_t\nabla_{y_i} (\frak \bau^L_i-\frak u_i)\right\|_{\tL(0,T;H_i)}\\
+\|\curl(\bau_0^L-u_0)\|_{\tL(0,T;H)}+\sum_{i=1}^n\|\scurl_{y_i}(\bau^L_i -u_i)\|_{\tL(0,T;H_i)}\leq c((\Delta t)^2+h_L^s).
\end{multline*}
\end{proposition}
\bproof
The proof is similar to that of Proposition \ref{prop:fullsemidiscerror}. We note that
\[
\|\partial_tq^L_{0,m+1/2}\|_{H}=\left\|{q^L_{0,m+1}-q^L_{0,m}\over\Delta t}\right\|_{H}\le \sup_{t\in (0,T)}\left\|{\partial q^L_0\over\partial t}\right\|_{H}\le ch_L^s
\]
due to \eqref{eq:estqLfull}. Similarly
\[
\|\partial_t\nabla_{y_i}\fq^L_{i,m+1/2}\|_{H_i}\le ch_L^s.
\]
We then get the conclusion. 
\eproof

\subsection{Sparse tensor product finite elements} \label{sec:sparsetensor}
To define the sparse tensor product finite element spaces, we employ the following orthogonal projection 
\beqas
P^{l0}:L^2(D)\to V^l,\ \ \ 
P^{l0}_\#:L^2(Y)\to V^l_\#,
%&& P^{l1}_\#:H^1_\#(Y)\to V^l_\#,\\
%&& P^{lc}_\#:H_\#(\curl,D)\to W^l_\#
\eeqas
with the convention $P^{-10}=0$, $P^{-10}_\#=0$. The detail spaces are defined as
% $P^{-11}_\#=0$ and $P^{-1c}_\#=0$. 
\[
{\cal V}^l=(P^{l0}-P^{(l-1)0})V^{l},\ \ {\cal V}^l_\#=(P^{l0}_\#-P^{(l-1)0}_\#)V^l.
\]
We note that
\[
V^l=\bigoplus_{0\le i\le l}{\cal V}^i\mbox{\ \ and\ \ }V^l_\#=\bigoplus_{0\le i\le l} {\cal V}^i_\#.
\]
Therefore the full tensor product spaces $\bar W_i^L$ and $\bar V_i^L$ can be written as
\[
\bar W_i^L=\left(\bigoplus_{0\le l_0,\ldots,l_{i-1}\le L}{\cal V}^{l_0}\otimes{\cal V}^{l_1}_\#\otimes\ldots\otimes{\cal V}^{l_{i-1}}_\#\right)\otimes W^L_\#,
\]
and 
\[
\bar V_i^L=\left(\bigoplus_{0\le l_0,\ldots,l_{i-1}\le L}{\cal V}^{l_0}\otimes{\cal V}^{l_1}_\#\otimes\ldots\otimes{\cal V}^{l_{i-1}}_\#\right)\otimes V^L_\#.
\]
We define the sparse tensor product finite element spaces as
\[
\hat W_i^L=\bigoplus_{l_0+\ldots+l_{i-1}\le L}{\cal V}^{l_0}\otimes{\cal V}^{l_1}_\#\otimes\ldots\otimes{\cal V}^{l_{i-1}}_\#\otimes W^{L-(l_0+\ldots+l_{i-1})}_\#;
\]
\[
\hat V_i^L=\bigoplus_{ l_0+\ldots+l_{i-1}\le L}{\cal V}^{l_0}\otimes{\cal V}^{l_1}_\#\otimes\ldots\otimes{\cal V}^{l_{i-1}}_\#\otimes V^{L-(l_0+\ldots+l_{i-1})}_\#,
\]
and
\[
\hat\bV^L=W^L\otimes\hat W_1^L\otimes\ldots\otimes\hat W_n^L\otimes \hat V_1^L\otimes\ldots\otimes\hat V_n^L.
\]
The  spatially semidiscrete sparse tensor product  finite element approximating problem is: Find $\hat\bu^L(t)\in \hat\bV^L$ such that :

 \begin{eqnarray}
&&\int_D\int_\bY \left[b(x,\by)\left({\partial^2\hat u^L_0\over \partial t^2}(t)+\sum_{i=1}^n\nabla_{y_i} {\partial^2\hat{\frak u}^{L}_i\over \partial t^2}(t)\right)\cdot\left(\hat v_0^L+\sum_{i=1}^n\nabla_{y_i} \hat{\frak v}^{L}_i\right)\right.\nonumber\\
&&\left.\qquad+a(x,\by)\left(\curl \hat u^L_0(t)+\sum_{i=1}^n\scurl_{y_i}\hat{ u}^{L}_i(t)\right)\cdot\left(\curl \hat v^L_0+\sum_{i=1}^n\scurl_{y_i}\hat{ v}^{L}_i\right)\right]d\by dx\nonumber\\
&&=\int_D f(x)\cdot \hat v^L_0(x)dx
\label{eq:sparsesemidiscprob}
\end{eqnarray}
for all $\hbv^L=(\hav_0^L,\hav_1^L,\ldots,\hav_n^L,\hafv_1^L,\ldots,\hafv_n^L)\in\hbV^L$. 
To find an error estimate for the sparse tensor product finite element approximation we note the following results
\begin{lemma}\label{lem:4.1} For $w\in \hat{\cal H}_i^s$,
\[
\inf_{w^L\in \hat W_i^L}\|w-w^L\|_{L^2(D\times Y_1\times\ldots\times Y_{i-1},H_\#(\curl,Y_i))}\le cL^{i/2}h_L^s\|w\|_{\hat{\cal H}_i^s};
\]
for $w\in \hat{\frak H}_i^s$,
\[
\inf_{w^L\in \hat V_i^L}\|w-w^L\|_{L^2(D\times Y_1\times\ldots\times Y_{i-1},H^1_\#(Y_i))}\le cL^{i/2}h_L^s\|w\|_{\hat{\frak H}_i^s}.
\]
\end{lemma}
The proof of these results follow from that for sparse tensor product approximation in \cite{BungartzGriebel} and \cite{HSelliptic}. Therefore, for $\bw\in \hat{\bH}^s$
\[
\inf_{\bw^L\in\hat\bV^L}\|\bw-\bw^L\|_{\bV}\le cL^{n/2}h_L^s\|\bw\|_{\hat{\bH}^s}.
\] 
We then have the following result.
\bp
\label{prop:sparsesemidiscerror}
Assume that condition (\ref{eq:Compatibility2}) and Assumption \ref{assum:regularityab} hold, $D$ is a Lipschitz polygonal domain, ${\rm div}f\in L^\infty(0,T;L^2(D))$ and $g_0, g_1$ belong to $H^s(\curl,D)$. If  $g_0^L$ and $g_1^L$ are chosen so that
\[
%\be
\|g_0^L-g_0\|_{V} \leq cL^{n/2}h_L^s \text{ and }     \|g_1^L-g_1\|_{H} \leq cL^{n/2}h_L^s, 
%\label{eq:fullfullydiscinitials}
%\ee
\]
where $s\in (0,1]$ is the constant in Proposition \ref{prop:regularityd2u0dt2},
then the solution of the spatially semidiscrete approximating problem \eqref{eq:sparsesemidiscprob} satisfies
\beqas
&&  \left\|\frac{\partial (\hau_0^L-u_0)}{\partial t}\right\|_{L^\infty(0,T;H)}+\sum_{i=1}^n\left\|\nabla_{y_i}\frac{\partial (\hafu^L_i-\fu_i)}{\partial t}\right\|_{L^\infty(0,T; H_i)}\\
&&\qquad\qquad+\|\scurl (\hau_0^L-u_0)\|_{L^\infty(0,T;H)}+\sum_{i=1}^n\|\scurl_{y_i}(\hau_i^L-u_i)\|_{L^\infty(0,T;H_i)}\le cL^{n/2}h_L^s. 
\eeqas
\epr
The proof of this proposition is identical to that of Proposition \ref{prop:fullsemidiscerror}. 

The  fully  discrete sparse tensor finite element product problem is:
For $m=1,...,M$ find $\hbu^L_m=(\hau_{0,m}^L,\hau_{1,m}^L,...,\hau_{n,m}^L, \hafu_{1,m}^L ,...,\hafu_{n,m}^L) \in \hbV^L$ such that 
\begin{eqnarray}
&&\int_D\int_\bY b(x,\by)\left(\partial^2_t\hau^L_{0,m}+\sum_{i=1}^n\nabla_{y_i} \partial^2_t\hafu^L_{i,m}\right)\cdot\left(\hav_0^L+\sum_{i=1}^n\nabla_{y_i} \hafv^L_i\right)\nonumber\\
&&\qquad+a(x,\by)\left(\curl \hau^L_{0 ,m , 1/4}+\sum_{i=1}^n\scurl_{y_i}\hau^L_{i ,m , 1/4}\right)\cdot\left(\curl \hav^L_0+\sum_{i=1}^n\scurl_{y_i}\hav^L_i\right)d\by dx\nonumber\\
&&\qquad=\int_D f_{m,1/4}(x)\cdot \hav^L_0(x)dx
\label{eq:sparsefullydiscprob}
\end{eqnarray}
for all $\hbv^L=(\hav_0^L, \hav_1^L,..., \hav^L_n, \hafv_1^L,...,\hafv_n^L) \in \hbV^L.$

For the fully discrete problem, we have the following result
\begin{proposition}
\label{prop:sparsefullydiscerror}
Assume that condition (\ref{eq:Compatibility2}) and Assumption \ref{assum:regularityab} hold, $D$ is a Lipschitz polygonal domain, ${\rm div}f\in L^\infty(0,T;L^2(D))$ and $g_0, g_1$ belong to $H^s(\curl,D)$ 
where $s\in (0,1]$ is the constant in Proposition \ref{prop:regularityd2u0dt2}.
If the initial value $\hat\bu_1^L$ is chosen so that
\begin{align*}
\|\partial_tp^L_{0,1/2}\|_{L^2(D)}&+\sum_{i=1}^n\|\partial_t\nabla_{y_i}\fp^L_{i,1/2}\|_{L^2(D\times\bY)}+\|\curl p^L_{0,1/2}\|_{L^2(D)}\\
&+\sum_{i=1}^n\|\scurl_{y_i}p^L_{i,1/2}\|_{L^2(D\times\bY)}\le c((\Delta t)^2+L^{n/2}h_L^s),
\end{align*} 
then 
\begin{align*}
\|\partial_t\hau_0^L-\partial_tu_0\|_{\tL(0,T;H)}&+\sum_{i=1}^n\|\nabla_{y_i} (\hafu^L_i-\fu_i)\|_{\tL(0,T;H_i)}+\|\curl(\hau_0^L-u_0)\|_{\tL(0,T;H)}\\
&+\sum_{i=1}^n\|\scurl_{y_i}(\hau^L_i -u_i)\|_{\tL(0,T;H_i)}\leq c((\Delta t)^2+L^{n/2}h_L^s).
\end{align*}
\end{proposition}
The proof is identical to that of Proposition \ref{prop:fullfullydiscerror}.

\section{Numerical correctors}
We construct numerical correctors in this section. For two scale problems, we derive an explicit error for the corrector in terms of the microscale $\ep$ and the FE meshsize. For general multiscale problems, as a homogenization error is not available, we derive a corrector without an error estimate. We first review some results for analytic correctors. 
\subsection{Analytic homogenization errors and correctors}
%
%\begin{proposition}\label{prop:2sregularhomerror}
%For two-scale problems, assume that
%$g_0=0$, $g_1 \in H^1(D)\bigcap V$, 
%$f \in H^1 (0,T; H)$, 
%$a(x,y) \in C(\overline{D}, C(\overline{Y})), 
%$u_0 \in L^\infty(0,T;H^1(\curl; D))$, ${\partial u_0\over\partial t}\in L^\infty(0,T;H^1(\curl,D))$, ${\partial^2u_0\over\partial t^2}\in L^\infty(0,T;H(\curl,D)\bigcap H^1(D))$, $N^r \in C^1(\bar{D}, C(\bar{Y}))^3$, $\scurl_y N^r \in  C^1(\bar{D}, C(\bar{Y})) $ and $\omega^r \in C^1(\bar{D}, C(\bar{Y}))$ for all $r=1,2,3$. There exists a constant $c$ that does not depend on $\varepsilon$ such that
%\[
%\left\|\frac{\partial u^\varepsilon}{\partial t}-\left[\frac{\partial u_0}{\partial t}+\nabla_y\frac{\partial \frak u_1}{\partial t}\left(\cdot, \cdot, \frac{\cdot}{\varepsilon} \right ) \right]\right\|_{L^{\infty}(0,T; H)}+\left\|\curl u^\varepsilon-\left[\curl u_0+ \curl_y u_1 \left(\cdot, \cdot, \frac{\cdot}{\varepsilon} \right )\right]\right\|_{L^{\infty}(0,T;H)} \leq c \varepsilon^{1/2}.
%\]

%\end{proposition}
%
%
For two scale problems, for conciseness of notations, we denote the solutions to cell problems $N_1^r$ and $w_1^r$ as $N^r$ and $w^r$. We have the following homogenization error for two scale problems. This result generalizes the well known $O(\ep^{1/2})$ homogenization error in \cite{BLP} and \cite{JKO} to the case where the solution $u_0$ of the homogenized equation possesses low regularity. We derive this error for two scale Maxwell wave equations, but the proof works verbatim for the two scale elliptic equations in \cite{BLP} and \cite{JKO}. The proof is lengthy and complicated so we refer to \cite{Mwt} for details.    
\begin{proposition}\label{prop:nonregularhomerror}
Assume that
$g_0=0$, $g_1\in H^1(D)\bigcap W$, $f \in H^1 (0,T; H)$, %$a(x,y) \in C(\overline{D}, C(\overline{Y}))$, 
$u_0$, ${\partial u_0\over\partial t}$ and ${\partial^2u_0\over\partial t^2}$ belong to $L^\infty(0,T;H^s(\curl; D))$ for $0<s\le 1$, $N^r \in C^1(\overline{D}, C(\overline{Y}))^3$, $\curl_y N^r \in  C^1(\overline{D}, C(\overline{Y})),$ $w^r \in C^1(\overline{D}, C(\overline{Y}))$ for all $r=1,2,3$. There exists a constant $c$ that does not depend on $\varepsilon$ such that
\[
\left\|\frac{\partial u^\varepsilon}{\partial t}-\left[\frac{\partial u_0}{\partial t}+\nabla_y\frac{\partial \frak u_1}{\partial t}\left(\cdot, \cdot, \frac{\cdot}{\varepsilon} \right ) \right]\right\|_{L^{\infty}(0,T; H)}+\left\|\curl u^\varepsilon-\left[\curl u_0+ \curl_y u_1 \left(\cdot, \cdot, \frac{\cdot}{\varepsilon} \right )\right]\right\|_{L^{\infty}(0,T;H)} \leq c \varepsilon^{\frac{s}{1+s}}.
\]
\end{proposition}
For the case of more than two scales, an explicit homogenization error is not available. However, we can deduce correctors when $\ep_{i-1}/\ep_i$ is an integer for all $i=2,\ldots,n$. 
We define the operator ${\cal U}_n^\ep$ as
\beqas
{\cal U}_n^\ep(\Phi)(x)
=
\int_{Y_1}\cdots\int_{Y_n}\Phi\Bigl(\ep_1\Bigl[{x\over\ep_1}\Bigr]
+
\ep_1 t_1,{\ep_2\over\ep_1}\Bigl[{\ep_1\over\ep_2}\Bigl\{{x\over\ep_1}\Bigr\}\Bigr]
+
{\ep_2\over\ep_1}t_2,\cdots,
\\
{\ep_n\over\ep_{n-1}}\Bigl[{\ep_{n-1}\over\ep_n}\Bigl\{{x\over\ep_{n-1}}\Bigr\}\Bigr]
+
{\ep_n\over\ep_{n-1}}t_n,\Bigl\{{x\over\ep_n}\Bigr\}\Bigr)dt_n\cdots dt_1
\eeqas
for all functions $\Phi\in L^1(D\times\bY)$. In the two scale case, we denote ${\cal U}_n^\ep$ by ${\cal U}^\ep$. We note the following property.
\begin{lemma}
For each function $\Phi\in L^1(D\times\bY)$ we have 
\be
\int_{D^{\ep_1}}{\cal U}^\ep_n(\Phi)dx=\int_D\int_\bY\Phi(x,\by)d\by dx,
\label{eq:calUint}
\ee
where $D^{\ep_1}$ is the $2\ep_1$ neighbourhood of $D$. 
\end{lemma}
We refer to \cite{CDG} for a proof. 
We have the following corrector result for multiscale problems.
\bp\label{prop:mscorrector}
Assume that $g_0=0$, $g_1\in W$ and $f\in H^1(0,T;H)$. We have
\beqas
\lim_{\ep\to 0}\left\|{\partial\ue\over\partial t}-{\cal U}_n^\ep\left({\partial u_0\over\partial t}+\sum_{i=1}^n\nabla_{y_i}{\partial \fu_i\over\partial t}\right)\right\|_{L^\infty(0,T;H)}+
\left\|\curl\ue-{\cal U}^\ep_n\left(\curl u_0+\sum_{i=1}^n\scurl_{y_i} u_i\right)\right\|_{L^\infty(0,T;H)}=0.
\eeqas
\epr
The proof of these corrector results can be found in \cite{Mwt}. 
\begin{remark}
Generally, the energy of a multiscale wave equation does not always converge to the energy of the homogenized wave equation when $g_0\ne 0$. We therefore restrict our consideration to the case where $g_0=0$. As shown in \cite{BFM}, the corrector of a general two scale wave equation involves the solution of another multiscale equation in the domain $D$. However, the scale interacting terms in \eqref{eq:msprob} always form a part of the corrector.
\end{remark}
\subsection{Numerical correctors for two-scale problems}
We now establish numerical correctors with an explicit error estimate for two scale problems. We first note the following result.
\begin{lemma}\label{lem:Uepu1est}
Assume that ${\partial u_0\over\partial t}\in L^\infty(0,T;H^s(D)^3)$ and $\curl u_0\in L^\infty(0,T;H^s(D)^3)$, $N^r\in C^1(\bar D,C^1_\#(\bar Y))^3$ and $w^r\in C^1(\bar D,C^1_\#(\bar Y))$, r=1,2,3. Then
\[
 \sup_{t\in [0,T]}\int_D\left|\scurl_yu_1\left(t,x,{x\over\ep}\right)-{\cal U}^\ep(\scurl_yu_1(t,\cdot,\cdot))(x)\right|^2dx\le c\ep^{2s}
\]
and 
\[
\sup_{t\in [0,T]}\int_D\left|{\partial\over\partial t}\nabla_y\fu_1\left(t,x,{x\over\ep}\right)-{\cal U}^\ep\left({\partial\over\partial t}\nabla_y\fu_1(t,\cdot,\cdot)\right)(x)\right|^2dx\le c\ep^{2s}.
\]
\end{lemma}
%\bproof
%The proof is similar to that for Lemma \ref{lem:Uepcurlu1}
%\eproof
The proof of this result is similar to that for the time independent case in Appendix B of \cite{Tiep1}, which utilizes the ideas of the proof of Lemma 5.5 in \cite{HSmultirandom}. We then have the following numerical corrector results. 
\begin{theorem} \label{2scaleconvergence}
Assume that condition \eqref{eq:Compatibility2} and Assumption \ref{assum:regularityab} hold, with $g_0=0$ and ${\rm div}f\in L^\infty(0,T;L^2(D))$, $D$ is a Lipschitz polygonal domain, and that $g_1^L$ is chosen so that $\|g_1^L-g_1\|_H\le ch_L^s$ where $s\in (0,1]$ is the constant in Proposition \ref{prop:regularityd2u0dt2}. Then for the solution of the semidiscrete problem \eqref{eq:fullsemidiscprob} using the full tensor product FEs, we have
\beqas
&&\left\|{\partial\ue\over\partial t}-\left({\partial\bar u_0^L\over\partial t}+{\cal U}^\ep\left({\partial\over\partial t}\nabla_{y}\bar\fu_1^L\right)\right)\right\|_{L^\infty(0,T;H)}\\
&&\qquad+\left\|\curl\ue-\left(\curl \bar u_0^L+{\cal U}^\ep\left(\scurl_{y}\bar u_1^L\right)\right)\right\|_{L^\infty(0,T;H)}\le c\left(h_L^s+\ep^{s\over s+1}\right).
\eeqas
For the semidiscrete problem \eqref{eq:sparsesemidiscprob} using the sparse tensor product FEs, if $\|g_1^L-g_1\|_H\le cL^{1/2}h_L^s$, we have
\beqas
&&\left\|{\partial\ue\over\partial t}-\left({\partial\hat u_0^L\over\partial t}+{\cal U}^\ep\left({\partial\over\partial t}\nabla_{y}\hat\fu_1^L\right)\right)\right\|_{L^\infty(0,T;H)}\\
&&\qquad+\left\|\curl\ue-\left(\curl \hat u_0^L+{\cal U}^\ep\left(\scurl_{y}\hat u_1^L\right)\right)\right\|_{L^\infty(0,T;H)}\le c\left(L^{1/2}h_L^s+\ep^{s\over s+1}\right).
\eeqas
\end{theorem}
\bproof
With the hypothesis of the theorem, from Propositions \ref{prop:regularityNomega} and \ref{prop:regularityd2u0dt2}, the conditions of Theorem  \ref{prop:nonregularhomerror} hold. We then have from \eqref{eq:calUint}
\beqas
\left\|{\cal U}^\ep\left({\partial\over\partial t}\nabla_{y}\fu_1(t)-{\partial\over\partial t}\nabla_{y}\bar\fu_1^L(t)\right)\right\|_H\le \left\|{\partial\over\partial t}\nabla_{y}\fu_1(t)-{\partial\over\partial t}\nabla_{y}\bar\fu_1^L(t)\right\|_{L^2(D \times Y)^3}
\eeqas
and 
\beqas
&&\left\|{\cal U}^\ep\left(\curl u_0+\scurl_{y}u_1-\curl \bar u_0^L-\scurl_{y}\bar u_1^L\right)\right\|_{H}\\
&&\qquad\qquad\le\left\|\curl u_0+\scurl_{y}u_1-\curl \bar u_0^L-\scurl_{y}\bar u_1^L\right\|_{L^2(D \times Y)^3} .
\eeqas
We note that
\beqas
&&\left\|{\partial\ue\over\partial t}-\left({\partial\bar u_0^L\over\partial t}+{\cal U}^\ep\left({\partial\over\partial t}\nabla_{y}\bar\fu_1^L\right)\right)\right\|_{L^\infty(0,T;H)}\\
&&\le\left\|{\partial\ue\over\partial t}-\left({\partial u_0\over\partial t}+{\partial\over\partial t}\nabla_{y}\fu_1(\cdot,\cdot,{\cdot\over\ep})\right)\right\|_{L^\infty(0,T;H)}\\
&&+ \left\|{\partial u_0\over\partial t}-{\partial \bar u_0^L\over\partial t}\right\|_{L^\infty(0,T;H)}+ \left\|{\partial\over\partial t}\nabla_{y}\fu_1(\cdot,\cdot,{\cdot\over\ep})-{\cal U}^\ep\left({\partial\over\partial t}\nabla_{y}\fu_1\right)\right\|_{L^\infty(0,T;H)}\\
&&+\left\|{\cal U}^\ep\left({\partial\over\partial t}\nabla_{y}\fu_1\right)-{\cal U}^\ep\left({\partial\over\partial t}\nabla_{y}\bar\fu_1^L\right)\right\|_{L^\infty(0,T;H)}
\eeqas
From Proposition \ref{prop:nonregularhomerror}, we have
\[
\left\|{\partial\ue\over\partial t}-\left({\partial u_0\over\partial t}+{\partial\over\partial t}\nabla_{y}\fu_1(\cdot,\cdot,{\cdot\over\ep})\right)\right\|_{L^\infty(0,T;H)} \le c \ep ^{\frac{s}{s+1}}.
\]
From Proposition \ref{prop:fullsemidiscerror}, we have
\[
\left\|{\partial u_0\over\partial t}-{\partial \bar u_0^L\over\partial t}\right\|_{L^\infty(0,T;H)} \le c h^s_L.
\]
From Lemma \ref{lem:Uepu1est}, we have
\[
\left\|{\partial\over\partial t}\nabla_{y}\fu_1(\cdot,\cdot,{\cdot\over\ep})-{\cal U}^\ep\left({\partial\over\partial t}\nabla_{y}\fu_1\right)\right\|_{L^\infty(0,T;H)} \leq c \ep^s.
\]
We note that
\begin{align*}
\left\|{\cal U}^\ep\left({\partial\over\partial t}\nabla_{y}\fu_1\right)-{\cal U}^\ep\left({\partial\over\partial t}\nabla_{y}\bar\fu_1^L\right)\right\|_{L^\infty(0,T;H)} &\leq \left\|{\partial\over\partial t}\nabla_{y}\fu_1-{\partial\over\partial t}\nabla_{y}\bar\fu_1^L\right\|_{L^\infty(0,T; L^2(D\times Y)^3)}\leq c h^s_L.
\end{align*}
Thus
\[
\left\|{\partial\ue\over\partial t}-\left({\partial\bar u_0^L\over\partial t}+{\cal U}^\ep\left({\partial\over\partial t}\nabla_{y}\bar\fu_1^L\right)\right)\right\|_{L^\infty(0,T;H)} \le c \ep ^{\frac{s}{s+1}}+c h^s_L+c \ep^s+c h^s_L \le c(h_L^s+\ep^{s\over s+1}).
\]
Similarly,
\beqas
&&\left\|\curl\ue-\left(\curl \bar u_0^L+{\cal U}^\ep\left(\scurl_{y}\bar u_1^L\right)\right)\right\|_{L^\infty(0,T;H)}\\ 
&&\qquad\qquad\le\|\curl\ue-\curl u_0-\curl_y u_1(\cdot,\cdot,{\cdot\over\ep})\|_{L^\infty(0,T;H)}+\left \| \curl u_0-\curl \bar u_0^L \right \|_{L^\infty(0,T;H)}\\
&&\qquad\qquad\qquad +\left\|\curl_y u_1(\cdot,\cdot,{\cdot\over\ep})- {\cal U}^\ep\left(\scurl_{y} u_1\right) \right \|_{L^\infty(0,T;H)}+\left \| {\cal U}^\ep\left(\scurl_{y} u_1\right)- {\cal U}^\ep\left(\scurl_{y}\bar u_1^L\right)\right \|_{L^\infty(0,T;H)}\\
&&\qquad\qquad \le c\ep^{s\over s+1}+ c h^s_L+c\ep ^s+ch^s_L \le c(h_L^s+\ep^{s\over s+1}).
\eeqas
We then have the desired estimate. 

The proof for the semidiscrete sparse tensor finite element solution is similar.  
\eproof

For fully discrete problems, we have the following results.
\begin{theorem}
Assume that condition \eqref{eq:Compatibility2} and Assumption \ref{assum:regularityab} hold, with $g_0=0$, $g_1 \in H^s(\curl,D)$ and ${\rm div}f\in L^\infty(0,T;L^2(D))$, $D$ is a Lipschitz polygonal domain ($s\in (0,1]$ is the constant in Proposition \ref{prop:regularityd2u0dt2}). For the  fully discrete full tensor product FE problem \eqref{eq:fullfullydiscprob}, assume that  $\bar\bu_1^L$ is chosen so that
\begin{multline*}
\left\|\partial_tp^L_{0,1/2}\right\|_{L^2(D)^3}+\left\|\partial_t\nabla_{y}\fp^L_{1,1/2}\right\|_{L^2(D\times Y)^3}+\left\|\curl p^L_{0,1/2}\right\|_{L^2(D)^3}+\left\|\scurl_{y}p^L_{1,1/2}\right\|_{L^2(D\times Y)^3}\\\le c((\Delta t)^2+h_L^s),
\end{multline*}
then 
\beqas
&&\Delta t\max_{0 \le m <M}\left\|{\partial\ue\over\partial t}(t_m)-\partial_t\bau_{0, m+1/2}^L-{\cal U}^\ep(\partial_t\nabla_{y} \frak \bau^L_{1,m+1/2})\right\|_{H}\\
&&\qquad\qquad+\left\|\curl\ue-\curl\bau_0^L-{\cal U}^\ep(\scurl_{y}\bau^L_1)\right\|_{\tL(0,T;H)}\leq c\left((\Delta t)^2+h_L^s+\ep^{s\over s+1}\right).
\eeqas
For the sparse tensor product FE problem \eqref{eq:sparsefullydiscprob}, if $\hat\bu_1^L$ is chosen so that
\begin{multline*}
\|\partial_tp^L_{0,1/2}\|_{L^2(D)^3}+\|\partial_t\nabla_{y}\fp^L_{1,1/2}\|_{L^2(D \times Y)^3}+\|\curl p^L_{0,1/2}\|_{L^2(D)^3}\\+\|\scurl_{y}p^L_{1,1/2}\|_{L^2(D \times Y)^3}\le c\left((\Delta t)^2+L^{1/2}h_L^s\right),
\end{multline*}
then 
\beqas
&&\Delta t\max_{0 \le m < M}\left\|{\partial\ue\over\partial t}(t_m)-\partial_t\hau_{0,m+1/2}^L-{\cal U}^\ep(\partial_t\nabla_{y} \hafu^L_{1,m+1/2})\right\|_{H}\\
&&\qquad+\left\|\curl\ue-\curl\hau_0^L-{\cal U}^\ep(\scurl_{y}\hau^L_1)\right\|_{\tL(0,T;H)}\leq c\left((\Delta t)^2+L^{1/2}h_L^s+\ep^{s\over s+1}\right).
\eeqas
\end{theorem}
\bproof
From the compatibility condition ${\partial ^4 u_0 \over \partial t^4} \in L^\infty (0,T;H)$ so ${\partial ^2 u_0 \over \partial t^2} \in C([0,T];H)$. To use the homogenization error in Theorem \ref{prop:nonregularhomerror}, we estimate
\[
{1\over\Delta t}(u_{0,m+1}-u_{0,m})-{\partial u_0\over\partial t}(t_m)={\partial u_0\over\partial t}(\tau)-{\partial u_0\over\partial t}(t_m)=\int_{t_m}^\tau{\partial^2u_0\over\partial t^2}(\sigma)d\sigma,
\]
for a value $t_m\le\tau\le t_{m+1}$. With the compatibility condition \eqref{eq:Compatibility2}, we have that ${\partial^2u_0\over\partial t^2}\in L^\infty(0,T;H)$. Thus
\[
\sup_{0\le m<M}\left\|\partial_tu_{0,m+1/2}-{\partial u_0\over\partial t}(t_m)\right\|_H\le c\Delta t.
\]
Similarly, using the smoothness of $N^r$ and $w^r$ for $r=1,2,3$, we have that ${\partial^2\over\partial t^2}\nabla\fu_1\in L^\infty(0,T;L^2(D\times Y))$. We note that
\begin{align*}
\partial_t\nabla_y\fu_{1,m+1/2}-{\partial\over\partial t}\nabla_{y}\fu_1(t_m)=&{\nabla_y\fu_{1,m+1}-\nabla_y\fu_{1,m}\over \Delta t}-{\partial\over\partial t}\nabla_{y}\fu_1(t_m)\\
=&{\partial\over\partial t}\nabla_{y}\fu_1(\tau)-{\partial\over\partial t}\nabla_{y}\fu_1(t_m) \\
=&\int_{t_m}^\tau{\partial^2\over\partial t^2}\nabla_{y}\fu_1(\sigma)d\sigma,
\end{align*}
for a value $t_m\le\tau\le t_{m+1}$. Thus
\[
\sup_{0\le m<M}\left\|\partial_t\nabla_y\fu_{1,m+1/2}-{\partial\over\partial t}\nabla_{y}\fu_1(t_m)\right\|_{L^2(D\times Y)}\le c\Delta t.
\]
We then get the result from Proposition \ref{prop:sparsefullydiscerror} and Theorem \ref{prop:nonregularhomerror}. \eproof
\subsection{Numerical correctors for multiscale problems}
As an explicit homogenization error is not available for the case of more than two scales, we do not distinguish the full and sparse tensor FE. We work with general FE spaces instead. For the semidiscrete problem \eqref{eq:nmsprob} we have:
\begin{theorem}
Assume that condition \eqref{eq:Compatibility1} holds with $g_0=0$, and $g_1^L$ is chosen so that $\lim_{L\to\infty}\|g_1^L-g_1\|_H=0$. Then the solution of problem \eqref{eq:nmsprob} satisfies
\begin{multline*}
\lim_{\substack{L\to\infty\\\ep \to 0}}\left\|{\partial\ue\over\partial t}-{\partial u_0^L\over\partial t}-{\cal U}^\ep_n\left(\sum_{i=1}^n{\partial\over\partial t}\nabla_{y_i}\fu_i^L\right)\right\|_{L^\infty(0,T;H)}\\+\left\|\curl\ue-\curl u_0^L-{\cal U}^\ep_n\left(\sum_{i=1}^n\curl_{y_i}u_i^L\right)\right\|_{L^\infty(0,T;H)}=0.
\end{multline*}
\end{theorem}
\bproof
The result is a direct consequence of Propositions \ref{prop:semidiscerror} and \ref{prop:mscorrector}.
Indeed,
\begin{align*}
&\left\|{\partial\ue\over\partial t}-{\partial u_0^L\over\partial t}-{\cal U}^\ep_n\left(\sum_{i=1}^n{\partial\over\partial t}\nabla_{y_i}\fu_i^L\right)\right\|_{L^\infty(0,T;H)}\\
&\quad\leq \left\|{\partial\ue\over\partial t}-{\partial u_0\over\partial t}-{\cal U}^\ep_n\left(\sum_{i=1}^n{\partial\over\partial t}\nabla_{y_i}\fu_i\right)\right\|_{L^\infty(0,T;H)}\\
&\qquad+\left \| {\partial u_0\over\partial t}-{\partial u_0^L\over\partial t} \right \|_{L^\infty(0,T;H)}+\left \|{\cal U}^\ep_n\left(\sum_{i=1}^n{\partial\over\partial t}\nabla_{y_i}\fu_i\right)-{\cal U}^\ep_n\left(\sum_{i=1}^n{\partial\over\partial t}\nabla_{y_i}\fu_i^L\right)  \right \|_{L^\infty(0,T;H)}.
\end{align*}
From Proposition \ref{prop:mscorrector}, we deduce that
\[
\left\|{\partial\ue\over\partial t}-{\partial u_0\over\partial t}-{\cal U}^\ep_n\left(\sum_{i=1}^n{\partial\over\partial t}\nabla_{y_i}\fu_i\right)\right\|_{L^\infty(0,T;H)} \to 0
\]
as $\ep \to 0$. From Proposition \ref{prop:semidiscerror} we deduce that $\left \| {\partial u_0\over\partial t}-{\partial u_0^L\over\partial t} \right \|_{L^\infty(0,T;H)} \to 0$ as $L\to \infty.$ The last term
\beqas
&&\left \|{\cal U}^\ep_n\left(\sum_{i=1}^n{\partial\over\partial t}\nabla_{y_i}\fu_i\right)-{\cal U}^\ep_n\left(\sum_{i=1}^n{\partial\over\partial t}\nabla_{y_i}\fu_i^L\right)  \right \|_{L^\infty(0,T;H)}
\le \left \|\sum_{i=1}^n{\partial\over\partial t}\nabla_{y_i}\fu_i-\sum_{i=1}^n{\partial\over\partial t}\nabla_{y_i}\fu_i^L  \right \|_{L^\infty(0,T;H_n)} \to 0\\
\eeqas
as $L \to \infty$. Thus
\[
\lim_{\substack{L\to\infty\\ \ep \to 0}}\left\|{\partial\ue\over\partial t}-{\partial u_0^L\over\partial t}-{\cal U}^\ep_n\left(\sum_{i=1}^n{\partial\over\partial t}\nabla_{y_i}\fu_i^L\right)\right\|_{L^\infty(0,T;H)}=0.
\]
Similarly, we have 
\begin{align*}
&\left \|\curl\ue-\curl u_0^L-{\cal U}^\ep_n\left(\sum_{i=1}^n\curl_{y_i}u_i^L\right)\right \|_{L^\infty(0,T;H)} \\
&\quad \le 
\left \|\curl\ue-\curl u_0-{\cal U}^\ep_n\left(\sum_{i=1}^n\curl_{y_i}u_i\right)\right \|_{L^\infty(0,T;H)}\!\!+\left \| \curl u_0-\curl u_0^L \right \|_{L^\infty(0,T; H)}\\
&\qquad+\left \| {\cal U}^\ep_n\left(\sum_{i=1}^n\curl_{y_i}u_i\right)-{\cal U}^\ep_n\left(\sum_{i=1}^n\curl_{y_i}u_i^L\right) \right \|_{L^\infty(0,T; H)}
\end{align*}
which tends to 0 as $L \to \infty$ and $\ep \to 0$. We then get the conclusion.
\eproof

For the fully discrete problem \eqref{eq:fullydiscprob} we have:

\begin{theorem} Assume that condition \eqref{eq:Compatibility1} holds with $g_0=0$, $\bu_1^L$ is chosen such that
\[
\lim_{L\to\infty}\|\partial_t p^L_{0,1/2}\|_{H}+\sum_{i=1}^n\|\partial_t\nabla_{y_i} \fp^L_{i,1/2}\|_{H_i}+\|\curl p^L_{0,1/2}\|_H+\sum_{i=1}^n\|\scurl_{y_i}p^L_{i,1/2}\|_{H_i}=0.
\]
Then 
\begin{multline*}
\lim_{\substack{L\to\infty\\\Delta t\to 0\\ \ep \to 0}}\sup_{0\le m<M}\left\|{\partial\ue\over\partial t}(t_m)-\partial_t u^L_{0,m+1/2}-{\cal U}^\ep_n\left(\sum_{i=1}^n\partial_t\nabla_{y_i}\fu^L_{i,m+1/2}\right)\right\|_H\\
+\left\|\curl\ue-\curl u^L_{0}-{\cal U}^\ep_n\left(\sum_{i=1}^n\scurl_{y_i}u^L_{i}\right)\right\|_{\tL(0,T;H_i)}=0.
\end{multline*}
\end{theorem}
\bproof
We have
\beqas
&&\left\|{\partial\ue\over\partial t}(t_m)-\partial_t u^L_{0,m+1/2}-{\cal U}^\ep_n\left(\sum_{i=1}^n\partial_t\nabla_{y_i}\fu^L_{i,m+1/2}\right)\right\|_H\\
&&\qquad\le\left \| {\partial\ue\over\partial t}(t_m)-{\partial u_0\over\partial t}(t_m)-{\cal U}_n^\ep\left(\sum_{i=1}^n\nabla_{y_i}{\partial \fu_i\over\partial t}(t_m)\right) \right\|_H+\left \|{\partial u_0\over\partial t}(t_m)-\partial_t u_{0,m+1/2}  \right \|_H\\
&&\qquad\qquad+\left \|{\cal U}_n^\ep\left(\sum_{i=1}^n\nabla_{y_i}{\partial \fu_i\over\partial t}(t_m)\right)-  {\cal U}^\ep_n\left(\sum_{i=1}^n\partial_t\nabla_{y_i}\fu_{i,m+1/2}\right)\right \|_H
+\left \|\partial_t u_{0,m+1/2}-\partial_t u^L_{0,m+1/2}  \right \|_H\\
&&\qquad\qquad+\left \| {\cal U}^\ep_n\left(\sum_{i=1}^n\partial_t\nabla_{y_i}\fu_{i,m+1/2}\right)-{\cal U}^\ep_n\left(\sum_{i=1}^n\partial_t\nabla_{y_i}\fu^L_{i,m+1/2}\right) \right \|_H.
\eeqas
From Proposition \ref{prop:mscorrector} we deduce that 
$$\left \| {\partial\ue\over\partial t}(t_m)-{\partial u_0\over\partial t}(t_m)-{\cal U}_n^\ep\left(\sum_{i=1}^n\nabla_{y_i}{\partial \fu_i\over\partial t}(t_m)\right) \right\|_H \to 0 $$
as $\ep \to 0$.
As ${\partial^2u_0\over\partial t^2}\in L^\infty(0,T;H)$,
\[
\lim_{\Delta t\to 0}
\sup_{0\le m<M}\left\|\partial_tu_{0,m+1/2}-{\partial u_0\over\partial t}(t_m)\right\|_H=0;
\]
and from \eqref{eq:fui} we have ${\partial^2\over\partial t^2}\nabla_{y_i}\fu_i\in L^\infty(0,T;H_i)$ so
\[
\lim_{\Delta t\to 0}\sup_{0\le m<M}\left\|\partial_t\nabla_{y_i}\fu_{i,m+1/2}-{\partial\over\partial t}\nabla_{y_i}\fu_i(t_m)\right\|_H=0.
\]
Thus
\beqas
&&\sup_{0 \le m < M}\left \|{\cal U}_n^\ep\left(\sum_{i=1}^n\nabla_{y_i}{\partial \fu_i\over\partial t}(t_m)\right)-  {\cal U}^\ep_n\left(\sum_{i=1}^n\partial_t\nabla_{y_i}\fu_{i,m+1/2}\right)\right \|_H \\
&&\qquad\le \sup_{0 \le m < M}\left \|\sum_{i=1}^n\left(\nabla_{y_i}{\partial \fu_i\over\partial t}(t_m)-\partial_t\nabla_{y_i}\fu_{i,m+1/2}
\right)  \right \|_{L^2(D \times \bY)^3} \to 0 
\eeqas
as $\Delta t \to 0$. From the FE convergence, we have
\[
\sup_{0 \le m < M}\left \|\partial_t u_{0,m+1/2}-\partial_t u^L_{0,m+1/2}  \right \|_H \to 0
\]
as $L \to \infty, \Delta t \to 0$, and
\beqas
&&\sup_{0 \le m < M}\left \| {\cal U}^\ep_n\left(\sum_{i=1}^n\partial_t\nabla_{y_i}\fu_{i,m+1/2}\right)-{\cal U}^\ep_n\left(\sum_{i=1}^n\partial_t\nabla_{y_i}\fu^L_{i,m+1/2}\right) \right \|_H\\
&&\qquad\le \sup_{0 \le m < M}\left \|\sum_{i=1}^n \left (\partial_t\nabla_{y_i}\fu_{i,m+1/2}- \nabla_{y_i}\fu^L_{i,m+1/2}\right)  \right \|_{L^2(D \times \bY)^3} \to 0
\eeqas
as $L \to \infty, \Delta t \to 0$. Thus
$$\lim_{\substack{L\to\infty\\\Delta t\to 0\\ \ep \to 0}}\sup_{0\le m<M}\left\|{\partial\ue\over\partial t}(t_m)-\partial_t u^L_{0,m+1/2}-{\cal U}^\ep_n\left(\sum_{i=1}^n\partial_t\nabla_{y_i}\fu^L_{i,m+1/2}\right)\right\|_H=0.$$
Using similar argument, we have
\[
\lim_{\substack{L\to\infty\\\Delta t\to 0\\ \ep \to 0}}
\left\|\curl\ue-\curl u^L_{0}-{\cal U}^\ep_n\left(\sum_{i=1}^n\scurl_{y_i}u^L_{i}\right)\right\|_{\tL(0,T;H)}=0.
\]
We then get the conclusion.
\eproof

\section{Numerical results}

We present in this section some numerical examples for two scale problems that confirm our analysis.

%The detail spaces ${\cal V}^l$ and ${\cal V}^l_\#$ are defined via orthogonal projection in Section \ref{sec:sparsetensor}, which are difficult to construct in numerical implementation. 
To identify the detailed spaces defined in Subsection \ref{sec:sparsetensor}, we employ Riesz basis and define the equivalent norms in the spaces $L^2(D)$ and $L^2(Y)$. The Riesz basis functions satisfy:
\begin{assumption}\label{FEassumption}
(i) For all vectors $j \in \mathbb{N}_0^d$, there exists an index set $I^j \subset \mathbb{N}^d_0$ and
a set of basis functions $\phi^{jk}\in L^2(D)$ for $k\in I^j$, such that $V^l = \text{span}\left\{\phi^{jk} : |j|_{\infty}\le l\right\}$. 
For all $\phi = \sum_{|j|_{\infty}\leq l,k\in I^j}\phi^{jk}c_{jk}\in V^l$
\[
c_1 \sum_{\substack{  |j|_{\infty}\leq l \\   k\in I^j  }}  |c_{jk}|^2 \leq \| \phi \|^2_{L^2(D)} \leq c_2\sum_{\substack{  |j|_{\infty}\leq l \\   k\in I^j  }} |c_{jk}|^2, 
\]
where $c_1>0$ and $c_2>0$ are independent of $\phi$ and $l$. 
%A similar assumption holds for $H^1_0(D)$ for basis functions $\phi^{jk}\in H^1_0(D)$. 

(ii) For the space $L^2(Y)$, for all $j \in \mathbb{N}_0^d$, there exists 
an index set  $I^j_0 \subset \mathbb{N}_0^d$ and 
a set of periodic basis  functions $\phi^{jk}_0\in L^2(Y)$, $k\in I^j_0$, such that $V^l_\# = \text{span}\{\phi^{jk}_0 : |j|_{\infty}\le l\}$. For all $\phi = \sum_{|j|_{\infty}\leq l,k\in I^j_0}\phi^{jk}_0c_{jk}\in V^l_\#$
\[
c_3 \sum_{\substack{  |j|_{\infty}\leq l \\   k\in I^j_0  }} |c_{jk}|^2 \leq \| \phi \|^2_{L^2(Y)} \leq c_4\sum_{\substack{  |j|_{\infty}\leq l \\   k\in I^j_0  }}  |c_{jk}|^2
\]
where $c_3>0$ and $c_4>0$ are independent of $\phi$ and $l$. 
\end{assumption}
With respect to the norm equivalence, we define the detailed spaces as ${\cal V}^l=\text{span}\{\phi^{jk} : |j|_{\infty}= l\}$ and ${\cal V}^l_\#=\text{span}\{\phi^{jk}_0 : |j|_{\infty}= l\}$.

%(iii) For periodic functions,  for each $j \in \mathbb{N}_0^d$, there exists 
%a set of indices $I^j_1 \subset \mathbb{N}_0^d$ and 
%a set of basis functions $\phi^{jk}_{1}\in H^1_{\#}(Y)/\IR$, $k\in I_{1}^j\subset \mathbb{N}_0^d$, such that $\{\phi^{jk}_1 : |j|_{\infty}\le l\}$ is a basis of $V^l_\#/\IR$. %There are constants $c_6>c_5>0$ such that
%If $\phi = \sum_{|j|_{\infty}\leq l,k\in I^j_1}\phi^{jk}_1c_{jk}$, then
%\[
%c_5 \sum_{\substack{  |j|_{\infty}\leq l \\   k\in I^j_{1}  }}  |c_{jk}|^2 \leq \| \phi \|^2_{H^1_{\#}(Y)/\IR} \leq c_6\sum_{\substack{  |j|_{\infty}\leq l \\   k\in I^j_{1}  }}  |c_{jk}|^2
%\]
%where $c_5$ and $c_6$ are independent of $\phi$ and $l$.
%\end{assumption}
 
%Multi-dimensional bases may be derived by takine product of one-dimensional bases.

{\bf Example}
(i)  For the space $L^2(0,1)$, a Riesz basis can be constructed as follows.  Level $0$ contains three piecewise linear basis functions: $\psi^{01}$ obtains values $(1,0)$ at $(0,1/2)$ and is 0 in $(1/2,1)$, $\psi^{02}$ obtains values $(0, 1, 0)$ at $(0, 1/2, 1)$, and $\psi^{03}$ obtains values $(0,1)$ at $(1/2, 1)$ and is 0 in $(0,1/2)$. For other levels, the basis functions are constructed from the function $\psi$ that takes values $(0,-1,2,-1,0)$ at $(0,1/2,1,3/2,2)$, the left boundary function $\psi^{left}$ taking values $(-2,2,-1,0)$ at $(0,1/2,1,3/2)$, and the right boundary function $\psi^{right}$ taking values $(0, -1,2,-2)$ at $(1/2,1,3/2,2)$. For levels $j\geq 1$ with $I^j=\{1,2,\ldots,2^j\}$, the basis functions are  $\psi^{j1}(x) = 2^{j/2}\psi^{left}(2^j x)$, $\psi^{jk}(x)=2^{j/2}\psi(2^j x - k + 3/2)$ for $k = 2, \cdots, 2^j-1$ and $\psi^{j2^j} = 2^{j/2}\psi^{right}(2^j x - 2^j+2)$. This basis satisfies Assumption \ref{FEassumption} (i).\\
%(ii) For space $H^1_0(D)$, a simple hierarchical basis function that satisfies Assumption \ref{FEassumption} can be constructed from  the hat function that is piecewise linear and  obtains the values $(0,1,0)$ at $(0,1/2,1)$.
%At level $j$, $I^j=\{1,2,\ldots,2^l\}$ and  $\phi^{jk}(y)=2^{-j/2}\phi(2^jy-i+1)$. 

(ii) For $Y = (0,1)$, a  periodic Riesz basis for $L^2(Y)$ can be constructed by modifying the basis in (i). Level 0 %we exclude $\psi^{01}$, $\psi^{03}$ and include 
contains the periodic piecewise linear function that takes values $(1,0,1)$ at $(0,1/2,1)$ respectively. At other levels, the  functions  $\psi^{left}$ and $\psi^{right}$ are replaced by the piecewise linear functions that take values $(0,2, -1, 0)$ at $(0,1/2,1,3/2)$ and values $(0, -1,2,0)$ at $(1/2,1, 3/2 ,2)$ respectively.

%Another hierarchical basis  can be constructed from the simple hat function that obtains the values $(0,1,0)$ at $(0,1/2,1)$. The wavelet functions $\phi^{jk}(y)=2^{-j/2}\phi(2^jy-i+1)$.

A Riesz basis for the space $L^2((0,1)^d)$ can be constructed by taking the tensor products of the basis functions in $(0,1)$ with an appropriate scaling, see \cite{GO95}. 
\begin{remark}
We note that the norm equivalences above are not necessary for the approximations in Lemma \ref{lem:4.1} to hold, as explained in \cite{Hoangmonotone} and \cite{Tiep1}.
\end{remark}

In the first example, we consider a two scale Maxwell wave equation in the two dimension domain $D=(0,1)^2$. 
$$b^\varepsilon \frac{\partial^2 { u}^\varepsilon}{\partial t^2}+ \curl (a^\varepsilon \curl { u}^\varepsilon)={ f}(t,x), ~~~\text{in } D$$
$$u^\varepsilon (t, \cdot) \times \nu=0, \text{ on } \partial D$$
 $$u^\varepsilon (0, x)=0$$
  $$u_t^\varepsilon (0, x)=0$$
The coefficients are
$$a(x,y)=\frac{1}{(1+x_1)(1+x_2)(1+\cos^2 2 \pi y_1)(1+\cos^2 2 \pi y_2)},   $$
and
$$b(x,y)=\frac{(1+x_1)(1+x_2)}{(1+\cos^2 2 \pi y_1)(1+\cos^2 2 \pi y_2)}.$$
The exact homogenized coefficients are
$$a^0(x)=\frac{4}{9(1+x_1)(1+x_2)}$$
and 
$$b^0(x)= \frac{\sqrt{2}(1+x_1)(1+x_2)}{3}.$$
We choose
$$f(x_1,x_2)=\begin{pmatrix}
2\sqrt{2}(1+x_1)(1+x_2)x_1x_2(1-x_2)t+\dfrac{4t^3}{9(1+x_2)^2}\\ 
2\sqrt{2}(1+x_1)(1+x_2)x_1x_2(1-x_1)t+\dfrac{4t^3}{9(1+x_1)^2}
\end{pmatrix}$$
so that the solution to the homogenized equation is
$${u}_0=\begin{pmatrix}
x_1x_2(1-x_2)t^3\\ 
x_1x_2(1-x_1)t^3
\end{pmatrix}.$$
From the relation \eqref{eq:ui},
we compute the solution $\curl_y u_1$ exactly as
\begin{align*}
\curl u_1%&=\curl u_0 \curl_y \chi=\left({a^0(x)\over a(x,y)}-1)\right)\curl u_0\\
=\left({4(1+\cos^2 2 \pi y_1)(1+\cos^2 2 \pi y_2) \over 9}-1\right)(x_2-x_1)t^3.
\end{align*}
\begin{figure}
	\centering
	\includegraphics[width=\textwidth]{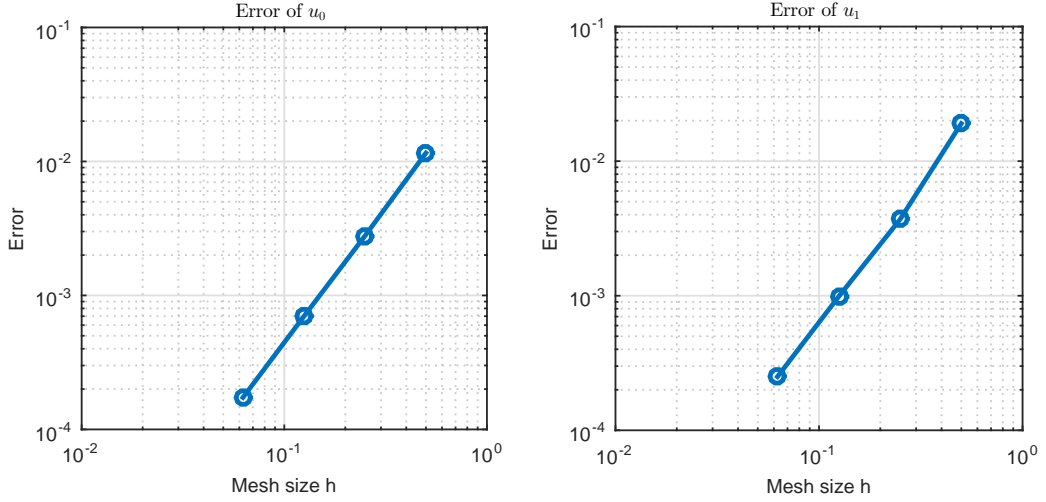}
	\caption{The sparse tensor  errors $\|u_0-u_0^L\|_{H(\curl,D)}$ and $\|\curl u_1-\curl u_1^L\|_{L^2(D)^3}$ }
\label{fig:5a}
\end{figure}
In Figure \ref{fig:5a} we plot the errors $\|u_0-u_0^L\|_{H(\curl,D)}$ and $\|\curl u_1-\curl u_1^L\|_{L^2(D)^3}$ versus the mesh size for the sparse tensor product FEs for  $(\Delta t, h)=(1/4, 1/4), (1/6, 1/8)$, $(1/8, 1/12)$ and $(1/16, 1/32)$. The result confirm  our analysis.

\begin{figure}
	\centering
	\includegraphics[width=\textwidth]{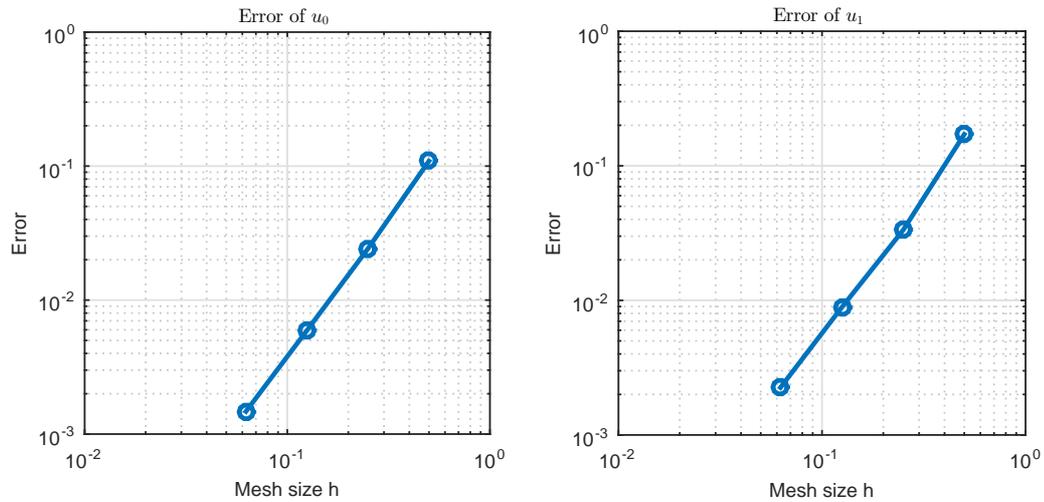}
	\caption{The sparse tensor  errors $\|u_0-u_0^L\|_{H(\curl,D)}$ and $\|\curl u_1-\curl u_1^L\|_{L^2(D)^3}$}
\label{fig:5b}
\end{figure}

In the second example, we choose
$$a(x,y)=\frac{(1+x_1)(1+x_2)}{(1+\cos^2 2 \pi y_1)(1+\cos^2 2 \pi y_2)},$$
and
$$b(x,y)=\frac{1}{(1+x_1)(1+x_2)(1+\cos^2 2 \pi y_1)(1+\cos^2 2 \pi y_2)}.$$
In this case, the homogenized coefficients are 
$$a^0(x)=\frac{4(1+x_1)(1+x_2)}{9}$$
and
$$b^0(x)=\frac{\sqrt{2}}{3(1+x_1)(1+x_2)}.$$
We choose
$$ f(x_1,x_2)=\begin{pmatrix}
\dfrac{2\sqrt{2}x_2(1-x_2)t}{(1+x_2)}+\dfrac{4t^3(1+x_1)(2x_2-x_1+1)}{3}\\ 
\dfrac{2\sqrt{2}x_1(1-x_1)t}{(1+x_1)}+\dfrac{4t^3(1+x_1)(2x_1-x_2+1)}{3}
\end{pmatrix}$$
so that the solution to the homogenized problem is
$$u_0=\begin{pmatrix}
(1+x_1)x_2(1-x_2)t^3\\ 
(1+x_2)x_1(1-x_1)t^3.
\end{pmatrix}$$
and
\begin{align*}
\curl u_1%&=\curl u_0 \curl_y \chi=\left({a^0(x)\over a(x,y)}-1)\right)\curl u_0\\
=\left({4(1+\cos^2 2 \pi y_1)(1+\cos^2 2 \pi y_2) \over 3}-1\right)(x_2-x_1)t^3.
\end{align*}
In Figure \ref{fig:5b} we plot the errors $\|u_0-u_0^L\|_{H(\curl,D)}$ and $\|\curl u_1-\curl u_1^L\|_{L^2(D)^3}$ versus the mesh size for the sparse tensor product FEs for  $(\Delta t, h)=(1/4, 1/4), (1/6, 1/8)$, $(1/8, 1/12)$ and $(1/16, 1/32)$. The result once more confirms our analysis.

{\bf Acknowledgement} The authors gratefully acknowledge a postgraduate scholarship of Nanyang Technological University, the  AcRF Tier 1 grant RG30/16, the Singapore A*Star SERC grant 122-PSF-0007 and the AcRF Tier 2 grant MOE 2013-T2-1-095 ARC 44/13. 

\bibliographystyle{plain}
\bibliography{references}
\end{document}